\RequirePackage[l2tabu, orthodox]{nag} 

\documentclass[reqno]{amsart}

\pdfoutput=1
\usepackage[utf8]{inputenc} 
\usepackage{bbm} 
\usepackage{commath} 
\usepackage{mathtools} 
\usepackage{todonotes} 
\usepackage[colorlinks, citecolor=red, breaklinks=true]{hyperref} 
\usepackage[stretch=10, shrink=10]{microtype} 
\usepackage[initials,msc-links,non-compressed-cites,nobysame]{amsrefs}[2007/10/22]

\newtheorem{theorem}{Theorem}[section]

\newtheorem{proposition}[theorem]{Proposition}
\newtheorem{fact}[theorem]{Fact}

\theoremstyle{definition}
\newtheorem{definition}[theorem]{Definition}

\theoremstyle{remark}
\newtheorem{remark}[theorem]{Remark}
\newtheorem{example}[theorem]{Example}

\numberwithin{equation}{section}

\newcommand{\ii}{\mathbbm{i}}
\newcommand{\jj}{\mathbbm{j}}
\newcommand{\kk}{\mathbbm{k}}

\DeclareMathOperator{\Span}{span}

\DeclareMathOperator{\cratio}{cr}
\DeclareMathOperator{\Id}{id}
\let\Im\relax
\DeclareMathOperator{\Im}{Im}

\usepackage{paralist,enumitem}
\setlist{listparindent=0pt,parsep=3pt}

\newcommand{\TitleWithUrl}[1]{\IfEmptyBibField{doi}%
  {\IfEmptyBibField{url}{\textit{#1}}%
    {\IfEmptyBibField{eprint}{\href {\BibField{url}}{\textit{#1}}}{\textit{#1}}}%
    }%
  {\href {https://doi.org/\BibField{doi}}{\textit{#1}}}}
\renewcommand{\eprint}[1]{\IfEmptyBibField{url}{\url{#1}}%
  {\href {\BibField{url}}{#1}}}

\BibSpec{article}{%
    +{}  {\PrintAuthors}                {author}
    +{,} { \TitleWithUrl}               {title}
    +{.} { }                            {part}
    +{:} { \textit}                     {subtitle}
    +{,} { \PrintContributions}         {contribution}
    +{.} { \PrintPartials}              {partial}
    +{,} { }                            {journal}
    +{}  { \textbf}                     {volume}
    +{}  { \PrintDatePV}                {date}
    +{,} { \issuetext}                  {number}
    +{,} { \eprintpages}                {pages}
    +{,} { }                            {status}
    +{,} { available at arXiv:\eprint}        {eprint}
    +{}  { \parenthesize}               {language}
    +{}  { \PrintTranslation}           {translation}
    +{;} { \PrintReprint}               {reprint}
    +{.} { }                            {note}
    +{.} {}                             {transition}
    +{}  {\SentenceSpace \PrintReviews} {review}
}

\BibSpec{collection.article}{%
    +{}  {\PrintAuthors}                {author}
    +{,} { \TitleWithUrl}                     {title}
    +{.} { }                            {part}
    +{:} { \textit}                     {subtitle}
    +{,} { \PrintContributions}         {contribution}
    +{,} { \PrintConference}            {conference}
    +{}  {\PrintBook}                   {book}
    +{,} { }                            {booktitle}
    +{,} { \PrintDateB}                 {date}
    +{,} { pp.~}                        {pages}
    +{,} { }                            {status}
    +{,} { available at \eprint}        {eprint}
    +{}  { \parenthesize}               {language}
    +{}  { \PrintTranslation}           {translation}
    +{;} { \PrintReprint}               {reprint}
    +{.} { }                            {note}
    +{.} {}                             {transition}
    +{}  {\SentenceSpace \PrintReviews} {review}
}
\BibSpec{book}{%
    +{}  {\PrintPrimary}                {transition}
    +{,} { \TitleWithUrl}               {title}
    +{.} { }                            {part}
    +{:} { \textit}                     {subtitle}
    +{,} { \PrintEdition}               {edition}
    +{}  { \PrintEditorsB}              {editor}
    +{,} { \PrintTranslatorsC}          {translator}
    +{,} { \PrintContributions}         {contribution}
    +{,} { }                            {series}
    +{,} { \voltext}                    {volume}
    +{,} { }                            {publisher}
    +{,} { }                            {organization}
    +{,} { }                            {address}
    +{,} { \PrintDateB}                 {date}
    +{,} { }                            {status}
    +{}  { \parenthesize}               {language}
    +{}  { \PrintTranslation}           {translation}
    +{;} { \PrintReprint}               {reprint}
    +{.} { }                            {note}
    +{.} {}                             {transition}
    +{}  {\SentenceSpace \PrintReviews} {review}
}

\BibSpec{thesis}{%
    +{}  {\PrintAuthors}                {author}
    +{.} { \TitleWithUrl}               {title}
    +{:} { \textit}                     {subtitle}
    +{,} { \PrintThesisType}            {type}
    +{,} { }                            {organization}
    +{} { \PrintDate}                  {date}
    +{,} { }                            {address}
    +{,} { \eprint}                     {eprint}
    +{,} { }                            {status}
    +{}  { \parenthesize}               {language}
    +{}  { \PrintTranslation}           {translation}
    +{;} { \PrintReprint}               {reprint}
    +{.} { }                            {note}
    +{.} {}                             {transition}
    +{}  {\SentenceSpace \PrintReviews} {review}
}

\title{Discrete constant mean curvature cylinders and isothermic tori}
\author{Joseph Cho}
\address[Joseph Cho]{Institute of Discrete Mathematics and Geometry, TU Wien, Wiedner Hauptstrasse 8-10/104, 1040 Wien, Austria}
\email{jcho@geometrie.tuwien.ac.at}

\author{Katrin Leschke}
\address[Katrin Leschke]{Department of Mathematics, University of Leicester, University Road, Leicester LE1 7RH, United
  Kingdom}
\email{k.leschke@leicester.ac.uk }

\author{Yuta Ogata}
\address[Yuta Ogata]{Department of Mathematics, Faculty of Science, Kyoto Sangyo University, 
  Motoyama, Kamigamo, Kita-ku, Kyoto-City, 603-8555, Japan.
 }
\email{yogata@cc.kyoto-su.ac.jp}


\begin{document}

\begin{abstract}
	We consider the monodromy problem of Darboux transforms of discrete isothermic surfaces using the integrable theory of discrete polarised curves.
	Then we provide, for the first time, closed-form discrete parametrisations of discrete isothermic cylinders, discrete constant mean curvature cylinders, and discrete isothermic tori.
\end{abstract}

\maketitle

\section{Introduction}

Discrete differential geometry aims to examine discrete surfaces through integrable systems technique (see, for example, \cite{bobenko_discrete_2008}), and has gained interest recently from both mathematical and applicational viewpoints.
However, a large number of results concerning structure-preserving discretisation of surfaces concern their local behaviour.
Here, we aim to consider a global aspect of discrete surfaces, focusing on the case of \emph{discrete isothermic surfaces}.


Isothermic surfaces \cite{bour_theorie_1862} constitute an integrable class of surfaces \cite{cieslinski_isothermic_1995} (see also \cite{burstall_curved_1997}), and admit spectral deformations \cite{bianchi_ricerche_1905, bianchi_complementi_1906, calapso_sulla_1903, calapso_sulle_1915} and \emph{Darboux transformations} \cite{darboux_sur_1899-1}.
In particular, Darboux transformations have been used to construct various isothermic surfaces with topological constraints by considering the \emph{monodromy problem}: notably, constant mean curvature (cmc) cylinders called bubbletons can be obtained as Darboux transformations of circular cylinders \cite{bianchi_lezioni_1903} (see also \cite{sterling_existence_1993}), while non-trivial isothermic tori were found via Darboux transforms of the homogeneous tori in $3$-sphere \cite{bernstein_non-special_2001}.

On the other hand, \emph{discrete isothermic surfaces} were defined in \cite{bobenko_discrete_1996-1} via the permutability of Darboux transformations \cite{bianchi_ricerche_1905}, and \emph{Darboux transformations} for discrete isothermic surfaces were defined in \cite{hertrich-jeromin_discrete_1999} using a discrete Riccati-type equation that is equivalent to a cross-ratios conditon.
A gauge theoretic approach to discrete isothermicity and its Darboux transforms was explored in both the quaternionic setup and the lightcone model \cite{hertrich-jeromin_transformations_2000, burstall_discrete_2014}, where isothermicity is characterised by a existence of a certain $1$-parameter family of discrete flat connections on the discrete trivial bundle, with Darboux transforms given by their parallel sections.
One can also consider the integrable reduction to the class of discrete cmc surfaces, so that certain Darboux transformations of a given discrete cmc surface result in another discrete cmc surfaces \cite{hertrich-jeromin_discrete_1999, burstall_discrete_2014}.

In this paper, we consider the monodromy problem for Darboux transformations of discrete isothermic surfaces to construct discrete isothermic cylinders, including discrete cmc cylinders, and discrete isothermic tori.
In particular, we are interested in obtaining \emph{closed-form discrete parametrisations} for explicit examples of such surfaces; thus, we will use the quaternionic setup for isothermicity, as this approach has proven to be a useful tool in obtaining explicit solutions to Riccati-type equations in both the smooth and discrete settings (see, for example, \cite{cho_new_2022, cho_periodic_nodate}).

In Section~\ref{sec:prelim}, we first review the gauge theoretic description of discrete isothermic surfaces and discrete polarised curves in the quaternionic setup; we will also briefly review the theory of discrete isothermic cmc surfaces.
Then in Section~\ref{sec:monodromy}, we view discrete isothermic surfaces as successive Darboux transformations of discrete polarised curves (a discrete analogue of the result obtained for semi-discrete isothermic surfaces in \cite{burstall_semi-discrete_2016}).
This change in viewpoint is crucial for solving the monodromy problem, as we show that the monodromy problem of Darboux transformations of discrete isothermic surfaces can be reduced to that of discrete polarised curves (see Theorem~\ref{thm:reduction}).

\begin{figure}
	\centering
	\begin{minipage}{0.47\linewidth}
		\centering
		\includegraphics[width=\textwidth]{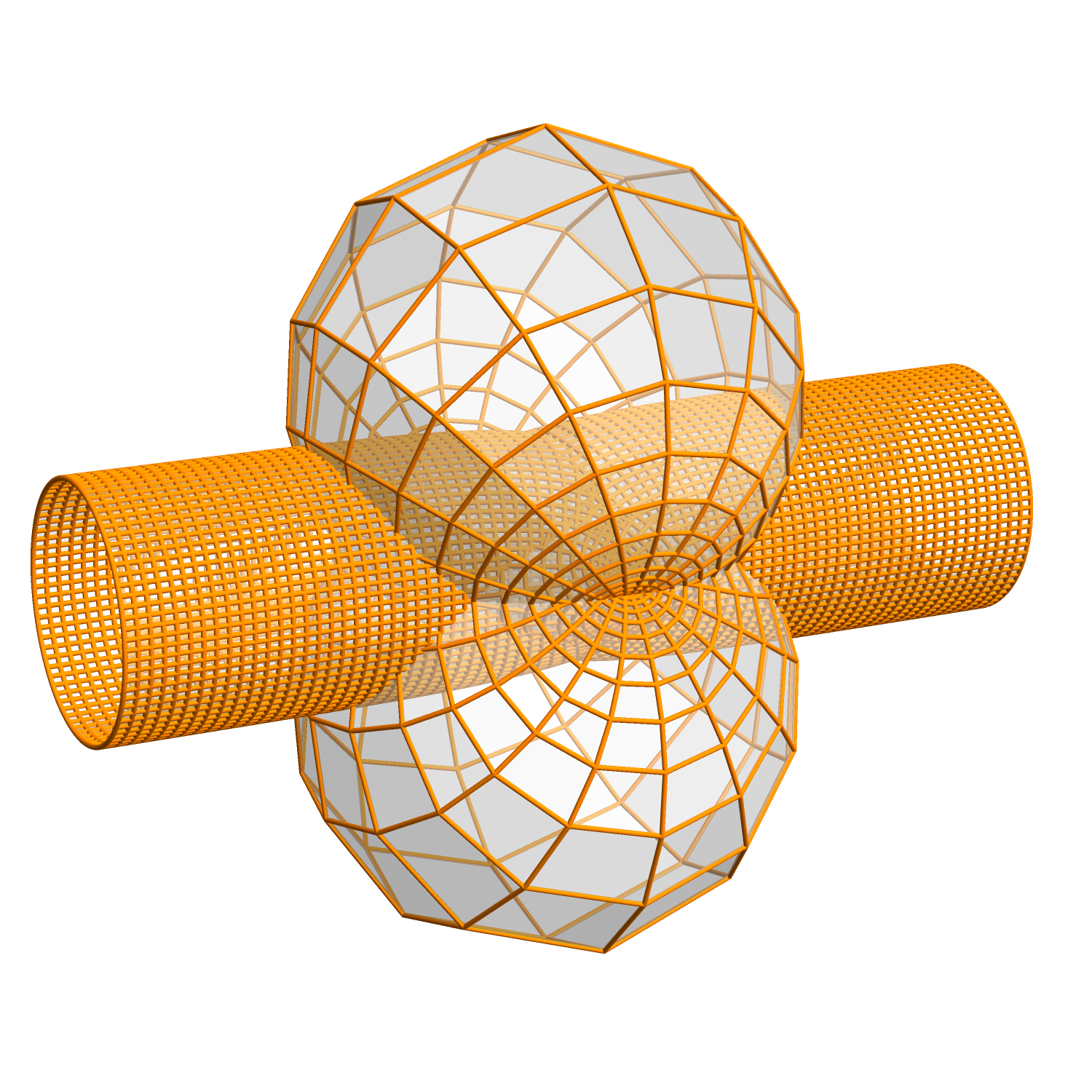}
	\end{minipage}
	\begin{minipage}{0.47\linewidth}
		\centering
		\includegraphics[width=\textwidth]{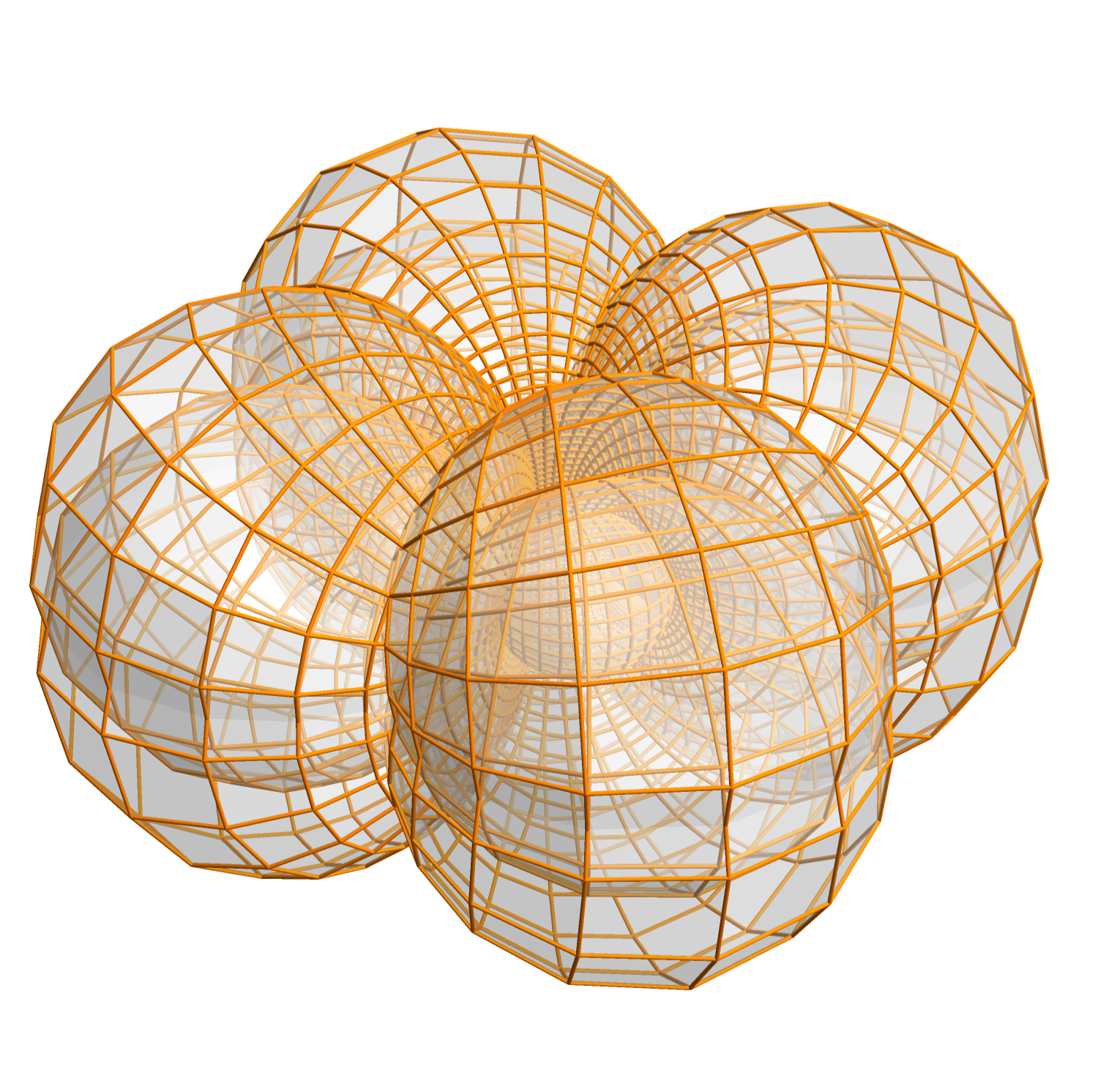}
	\end{minipage}
	\caption{Discrete constant mean curvature cylinder (on the left) and discrete isothermic torus (on the right) obtained as Darboux transformations of discrete isothermic surfaces.}
	\label{fig:introExam}
\end{figure}

In Section~\ref{sec:examples}, we consider Darboux transformations of explicit discrete isothermic surfaces.
After obtaining discrete isothermic cylinders as Darboux transformations of the discrete circular cylinder in Section~\ref{sect:circCyliso}, we use the fact that discrete circular cylinders are discrete cmc surfaces to find the Darboux transforms that are discrete cmc cylinders in Section~\ref{sect:circCylcmc}; these are the discrete analogues of the bubbletons.
Then we switch our attention to the Darboux transformations of homogeneous tori in $3$-sphere, and obtain explicit examples of discrete isothermic tori in Section~\ref{sect:tori}. (See also Figure~\ref{fig:introExam}.)

Importantly, we provide, for the first time, \emph{closed-form discrete parametrisations} of all the explicit examples we consider by solving the discrete Riccati-type equation (see \eqref{eqn:yay}, \eqref{eqn:expPar}, or \eqref{eqn:abinvtori}); the explicitness of our formulation allows us to observe in Section~\ref{sec:continuum} that under the appropriate continuum limit, the closure of Darboux transforms is preserved.
Thus, our work suggests that structure preserving discretisation can provide valuable insights in search of new examples of compact isothermic surfaces.

{\bf Acknowledgements.} The second author gratefully acknowledges the partial support of this work provided by the Research Institute for Mathematical Sciences, an International Joint Usage/Research Center located at Kyoto University.
The third author gratefully acknowledges the support from Grants-in-Aid of JSPS Research Fellowships for Young Scientist 21K13799. 

\section{Preliminaries}\label{sec:prelim}
We review the quaternionic vector bundle formalism of Darboux transformations of discrete isothermic surfaces and discrete polarised curves.

\subsection{Discrete flat connections}
\begin{figure}
		\centering
		\begin{tikzpicture}[scale=0.8]
			\node[circle, fill=black, inner sep = 0pt, outer sep = 0pt, minimum size = 5pt, label={[shift={(0,-0.8)}, font=\small]$(m,n)$}] (a) at (0,0) {};
			\node[circle, fill=black, inner sep = 0pt, outer sep = 0pt, minimum size = 5pt, label={[shift={(-1,-0.4)}, font=\small]$(m, n+1)$}] (b) at (0,2) {};
			\node[circle, fill=black, inner sep = 0pt, outer sep = 0pt, minimum size = 5pt, label={[shift={(0,-0.8)}, font=\small]$(m+1, n)$}] (c) at (2,0) {};
			\node[circle, fill=black, inner sep = 0pt, outer sep = 0pt, minimum size = 5pt, label={[shift={(-1,-0.4)}, font=\small]$(m, n+2)$}] (d) at (0,4) {};
			\node[circle, fill=black, inner sep = 0pt, outer sep = 0pt, minimum size = 5pt, label={[shift={(0,0)}, font=\small]$(m+1,n+1)$}] (e) at (2,2) {};
			\node[circle, fill=black, inner sep = 0pt, outer sep = 0pt, minimum size = 5pt, label={[shift={(0,-0.8)}, font=\small]$(m+2, n)$}] (f) at (4,0) {};
			\path[-, line width = 1pt]
				(a) edge (b)
				(a) edge (c)
				(b) edge (e)
				(c) edge (e);
			\path[line width=1pt, line cap=round, dash pattern=on 2pt off 4\pgflinewidth]
				(b) edge (d)
				(c) edge (f);
			\node[circle, fill=black, inner sep = 0pt, outer sep = 0pt, minimum size = 5pt, label={[shift={(-0.25,-0.6)}, font=\small]$i$}] (g) at (7,0) {};
			\node[circle, fill=black, inner sep = 0pt, outer sep = 0pt, minimum size = 5pt, label={[shift={(-0.25,-0.1)}, font=\small]$\ell$}] (h) at (7,2) {};
			\node[circle, fill=black, inner sep = 0pt, outer sep = 0pt, minimum size = 5pt, label={[shift={(0.2,-0.7)}, font=\small]$j$}] (i) at (9,0) {};
			\node[circle, fill=black, inner sep = 0pt, outer sep = 0pt, minimum size = 5pt, label={[shift={(-0.25,-0.3)}, font=\small]$p$}] (j) at (7,4) {};
			\node[circle, fill=black, inner sep = 0pt, outer sep = 0pt, minimum size = 5pt, label={[shift={(0.2,-0.1)}, font=\small]$k$}] (k) at (9,2) {};
			\node[circle, fill=black, inner sep = 0pt, outer sep = 0pt, minimum size = 5pt, label={[shift={(0,-0.7)}, font=\small]$q$}] (l) at (11,0) {};
			\node (m) at (12,0) {$\phantom{q}$};
			\path[-, line width = 1pt]
				(g) edge (h)
				(g) edge (i)
				(h) edge (k)
				(i) edge (k);
			\path[line width=1pt, line cap=round, dash pattern=on 2pt off 4\pgflinewidth]
				(h) edge (j)
				(i) edge (l);
		\end{tikzpicture}
		\caption{An elementary quadrilateral with auxiliary points and their labels.}
		\label{fig:elem}
	\end{figure}
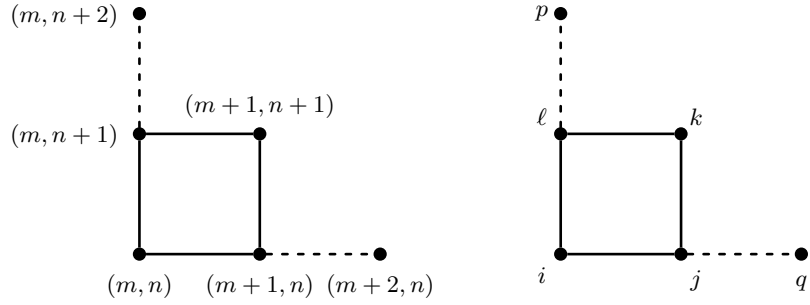

First, we briefly review the vector bundle theory for discrete systems (see, for example, \cite{burstall_notes_2017, burstall_discrete_2023}).
Let $(m,n) \in \Sigma^2 \subset \mathbb{Z}^2$ be a simply connected discrete domain in the sense of \cite{burstall_discrete_2023}.
For the labeling of vertices used throughout the paper, please refer to Figure~\ref{fig:elem}.
Writing $f(m,n) = f_{m,n}$ or $f(i) = f_i$ for any function $f$ defined on $\Sigma^2$, consider a discrete bundle $W$, that is, an assignment of sets $V_i$ to each vertex $i \in \Sigma^2$. We call $\sigma : \Sigma^2 \to \cup V_i$ a section, i.e.\ $\sigma \in \Gamma W$, if we have $\sigma_i \in V_i$ for all $i$.

A discrete connection is an assignment of bijections $r_{ji} : V_i \to V_j$ to each oriented edge $(ij)$ so that $r_{ji}r_{ij} = \Id_j$; discrete gauge transformations $G$ act on the discrete connections $r_{ji}$ via
	\[
		(G \bullet r)_{ji} = G_j r_{ji} G_i^{-1}
	\]
for automorphisms $G_i : V_i \to V_i$ defined on every vertex $i \in \Sigma$.

A connection is \emph{flat} if, on every elementary quadrilateral $(ijk\ell)$, we have 
	\[
		r_{i\ell}r_{\ell k}r_{kj}r_{ji} = \Id_i,
	\]
and it follows that any gauge transform of a flat connection is flat.
For a flat connection, we call a section $r$-parallel if $\sigma_{j} = r_{ji} \sigma_i$ on every edge $(ij)$.

\subsection{Integrable structure of discrete isothermic surfaces}
In this section. we review the basic integrable structure of discrete isothermic surfaces in $\mathbb{H} \cong \mathbb{R}^4$, thoroughly investigated in \cite{hertrich-jeromin_transformations_2000}.
While reviewing, we will adapt the results of \cite{hertrich-jeromin_transformations_2000} in terms of discrete vector bundle theoretic language.

\subsubsection{Discrete isothermic surfaces}
For a simply connected $\Sigma^2 \subset \mathbb{Z}^2$ with $(m,n) \in \Sigma^2$, we first recall the definition of discrete isothermic surfaces in the set of quaternions $\mathbb{H} := \Span_\mathbb{R}\{1, \ii, \jj, \kk\} \cong \mathbb{R}^4$ given in \cite[Definition 3.2]{hertrich-jeromin_transformations_2000} (see also \cite[Definition~6]{bobenko_discrete_1996-1}):

\begin{definition}\label{def:discIsothermic}
	A discrete surface $f : \Sigma^2 \to \mathbb{H}$ is called \emph{discrete isothermic} if
	\begin{equation}\label{eqn:isoCr}
		\cratio(f_i, f_j, f_k, f_\ell) := (f_i - f_j)(f_j - f_k)^{-1}(f_k - f_\ell)(f_\ell - f_i)^{-1} = \frac{\mu_{i\ell}}{\mu_{ij}}
	\end{equation}
	on every elementary quadrilateral $(ijk\ell)$ for some function $\mu$ defined on unoriented edges satisfying the edge-labeling property, i.e.\, $\mu_{ij} = \mu_{k\ell}$ and $\mu_{i\ell} = \mu_{jk}$.
    The function $\mu$ is called the \emph{cross-ratios factorising function}.
\end{definition}

\begin{remark}
	The definition of isothermicity implies that the image of every elementary quadrilateral $(ijk\ell)$ is \emph{concircular}.
	Therefore, the class of discrete isothermic surfaces is a special class of \emph{circular nets}, a discrete analogue of curvature line parametrised surfaces \cite[p.\ 220]{nutbourne_differential_1988}.
	For this reason, we will call $f_{m,n_0}$ or $f_{m_0, n}$ for some fixed $m_0$ or $n_0$ a \emph{discrete $m$-curvature line} or \emph{discrete $n$-curvature line}, respectively (see also, for example, \cite[Definition~2.7]{cho_discrete_2021}).
\end{remark}

For the discrete $1$-form $\omega_{ij}$ defined as
	\[
		\omega_{ij} = \frac{1}{\mu_{ij}}\dif{f}_{ij}^{-1}
	\]
on every edge $(ij)$ where
	\[
		\dif{f}_{ij} := f_i - f_j,
	\]
we have from \cite[Definition 3.4]{hertrich-jeromin_transformations_2000} (or \cite[Theorem~6]{bobenko_discrete_1996-1}) that $f$ is isothermic if and only if $\omega$ is closed, i.e.\ on every elementary quadrilateral $(ijk\ell)$,
	\[
		\omega_{ij} + \omega_{jk} + \omega_{k\ell} + \omega_{\ell i} = 0.
	\]
In this case, the dual surface $f^*$ of $f$, also called the \emph{Christoffel transform}, can be recursively defined via
	\begin{equation}\label{eqn:dual}
		\dif{f}_{ij}^* = \omega_{ij},
	\end{equation}
and is also isothermic with identical cross-ratios factorising function (\cite[p.\ 197]{bobenko_discrete_1996-1}) 
	\begin{equation}\label{eqn:isoCr2}
		\cratio(f^*_i, f^*_j, f^*_k, f^*_\ell) = \frac{\mu_{i\ell}}{\mu_{ij}}.
	\end{equation}
We note here that the Christoffel transform of a discrete isothermic surface is well-defined up to translations and homotheties.
This freedom can be observed in the choice of the cross-ratios factorising function $\mu$, as it can be globally scaled by a constant.

Now to introduce the vector bundle approach to discrete isothermic surfaces in the realm of M\"obius geometry, take the trivial bundles $\underline{\mathbb{H}}^2 := \Sigma^2 \times \mathbb{H}^2$ and $\underline{\mathbb{HP}}^1 := \Sigma^2 \times \mathbb{HP}^1$ where $\mathbb{H}^2$ is considered as a right quaternionic vector space.
Then the points $f \in \mathbb{H}$ can be identified with its lift $F \in \mathbb{HP}^1$ via
	\[
		\mathbb{H} \ni f \sim \begin{pmatrix}f \\ 1 \end{pmatrix} \mathbb{H} =: \psi \mathbb{H} =: F \in \mathbb{HP}^1.
	\]
With this identification, the lift $F : \Sigma^2 \to \mathbb{HP}^1$ of a discrete surface $f$ can be viewed as a quaternionic line subbundle of the trivial bundle $\underline{\mathbb{H}}^2$.

Now, for $D(\lambda)_{ji} : \{i\} \times \mathbb{H}^2 \to  \{j\} \times \mathbb{H}^2$ with $\lambda \in \mathbb{R}$ defined on every edge $(ij)$ via
	\begin{equation}\label{eqn:flat1}
		D(\lambda)_{ji} := \Id_{ji} + \begin{pmatrix} 0 & \dif{f}_{ij} \\ \lambda \omega_{ij} & 0 \end{pmatrix},
	\end{equation}
let $r(\lambda)_{ji} : \{i\} \times \mathbb{H}^2 \to  \{j\} \times \mathbb{H}^2$ be defined as
	\begin{equation}\label{eqn:flat2}
		r(\lambda)_{ji} := (G \bullet D(\lambda))_{ji} = G_j D(\lambda)_{ji} G_i^{-1}
	\end{equation}
where $G = (e \,\, \psi)$ for $e = \begin{psmallmatrix} 1 \\ 0 \end{psmallmatrix}$.
Then it is known that $r(\lambda)$ can be written as
	\begin{equation} \label{eqn:flatRep}
			r(\lambda)_{ji} = \Id_{ji} + \lambda \eta_{ji}
				:= \Id_{ji} + \lambda
					\begin{pmatrix}
						f_j \omega_{ij} & -f_j \omega_{ij} f_i \\
						\omega_{ij} & -\omega_{ij}f_i
					\end{pmatrix}.
	\end{equation}

The projective transformation $r(\lambda)^P_{ji} : \{i\} \times \mathbb{HP}^1 \to  \{j\} \times \mathbb{HP}^1$ induced by $r(\lambda)_{ji}$ defines a connection on the trivial bundle $\underline{\mathbb{HP}}^1$.
For an isothermic surface, connections in this $1$-parameter family of connections are flat:

\begin{fact}[{\cite[Equation (13)]{hertrich-jeromin_transformations_2000}}, \cite{cho_periodic_nodate}]\label{thm:flat}
	Let $f : \Sigma^2 \to \mathbb{H}$ be a discrete surface, with its lift $F : \Sigma^2 \to \mathbb{HP}^1$.
	Then $f$ is a discrete isothermic surface $f : \Sigma^2 \to \mathbb{H}$ if and only if $r(\lambda)^P$ is flat.
	Moreover, if $\mu$ is the cross-ratios factorising function of a discrete isothermic surface $f$, then $r(\lambda)$ can be written as
		\[
			r(\lambda)_{ji} = \pi_i + \frac{\mu_{ij}-\lambda}{\mu_{ij}} \pi_j
		\]
	on every edge $(ij)$ where $\pi_i$ denotes the projection onto $F_i$.
\end{fact}

\begin{remark}
    The $\eta_{ji}$ is called the \emph{retraction form} in \cite[Definition 3.6]{hertrich-jeromin_transformations_2000}.
\end{remark}

\subsubsection{Darboux transformations} On the other hand, the definition of Darboux transforms $\hat{f}$ of $f$ with spectral parameter $\nu$ is given in \cite{hertrich-jeromin_discrete_1999} via a cross-ratios condition:
\begin{definition}\label{def:darb}
	Let $f : \Sigma^2 \to \mathbb{H}$ be a discrete isothermic surface with cross-ratios factorising function $\mu$.
	A second discrete surface $\hat{f} : \Sigma^2 \to \mathbb{H}$ is called a \emph{Darboux transform of $f$ with spectral parameter $\nu$} if
		\[
			\cratio(f_i, f_j, \hat{f}_j, \hat{f}_i) = \frac{\nu}{\mu_{ij}}
		\]
	on every edge $(ij)$. The two surfaces $f$ and $\hat{f}$ are called a \emph{Darboux pair}, and share the same cross-ratios factorising function.
\end{definition}

A characterisation of Darboux transforms in terms of flat connections is given in \cite[Lemma~3.17]{hertrich-jeromin_transformations_2000}:
\begin{fact}
	For two isothermic surfaces $f, \hat{f} : \Sigma^2 \to \mathbb{H}$, let $F, \hat{F} : (\Sigma, \frac{1}{\mu}) \to \mathbb{HP}^1$ be the corresponding lifts.
	Denoting by $r(\lambda)^P$, the associated family of flat connections of $F$, we have that $\hat{F}$ is a Darboux transform of $F$ with spectral parameter $\mu$ if and only if $\hat{F}$ is $r(\nu)^P$--parallel, that is,
		\begin{equation}\label{eqn:parallel}
			r(\nu)^P_{ji} \hat{F}_i = \hat{F}_j
		\end{equation}
	on every edge $(ij)$.
\end{fact}

We remark here that the Darboux transforms of a discrete isothermic surface are uniquely determined by the choice of a spectral parameter and an initial condition.
Finally, we recall the condition for the Darboux transform $\hat{f}$ of a discrete isothermic surface $f$ in $\mathbb{S}^3$ to take values again in $\mathbb{S}^3$:
\begin{fact}[{\cite[p.\ 479]{hertrich-jeromin_transformations_2000}, (cf.\ \cite[Lemma~3.8]{cho_periodic_nodate})}] \label{fact:s3}
	Let $f : \Sigma^2 \to \mathbb{S}^3$ be a discrete isothermic surface with its Darboux transform $\hat{f}$.
	Then $\hat{f}$ takes values in $\mathbb{S}^3$ if and only if $\hat{f}_{m_0, n_0} \in \mathbb{S}^3$ for some $(m_0,n_0) \in \Sigma^2$.
\end{fact}

\subsubsection{Constant mean curvature surfaces}
An important subclass of isothermic surfaces is the class of constant mean curvature (cmc) surfaces in $\mathbb{R}^3$.
Discrete isothermic cmc surfaces have been defined using various characterisations in the smooth case.
Among many, we highlight a few that are particularly relevant to us:
\begin{itemize}
	\item surfaces with constant mean curvature defined via Steiner's formula \cite{schief_unification_2003, schief_maximum_2006, pottmann_geometry_2007, bobenko_curvature_2010},
	\item surfaces with Christoffel transform as parallel surfaces \cite{bobenko_discretization_1999, hertrich-jeromin_discrete_1999},
	\item surfaces admitting a simultaneous Christoffel and Darboux transform \cite{hertrich-jeromin_discrete_1999},
	\item surfaces admitting linear conserved quantities in the Minkowski $5$-space \cite{burstall_discrete_2014}.
\end{itemize}
These various definitions have been shown to be equivalent (see, for example, \cite{bobenko_discrete_2008, burstall_discrete_2014}).
For the purposes of the examination of Darboux transforms of discrete cmc surfaces, we will use the following definition:
\begin{definition}[{\cite[\S~5]{hertrich-jeromin_discrete_1999}}]\label{def:cmcDef}
	Let $f : \Sigma^2 \to \mathbb{R}^3 \cong \Im \mathbb{H}$ be a discrete isothermic surface.
	Then $f$ is a discrete cmc $H$ surface if it admits a Christoffel transform $f^* : \Sigma^2 \to \mathbb{R}^3$ so that
		\[
			|f^*_i - f_i|^2 = \frac{1}{H^2}
		\]
	at every vertex $i \in \Sigma^2$.
	In this case, the Christoffel transform $f^*$ is referred to as the parallel cmc surface of $f$.
\end{definition}
This definition of discrete cmc surfaces via parallel cmc surfaces allows us to recover another characterisation as follows:
\begin{fact}[{\cite[\S~5]{hertrich-jeromin_discrete_1999}}]\label{fact:cmcDarb}
	Let $f : \Sigma^2 \to \mathbb{R}^3 \cong \Im \mathbb{H}$ be a discrete cmc surface with parallel cmc surface $f^*$.
	Then $f^*$ is a Darboux transform of $f$ with spectral parameter $H^2$.
\end{fact}

We have seen that discrete isothermic cmc surfaces admit a $4$-parameter family of Darboux transforms that stay in $\mathbb{R}^3$.
Subject to a condition on the initial point, the Darboux transforms are again cmc surfaces:
\begin{fact}[{\cite[\S~5]{hertrich-jeromin_discrete_1999}}]\label{fact:cmcCond}
	Let $f : \Sigma^2 \to \mathbb{R}^3 \cong \Im \mathbb{H}$ be a discrete cmc-$H$ surface, and let $\hat{f}: \Sigma^2 \to \mathbb{R}^3 \cong \Im \mathbb{H}$ be its Darboux transform with spectral parameter $\nu$.
	Then $\hat{f}$ is a discrete cmc surface if and only if for some vertex $i \in \Sigma^2$, we have
		\[
			| \hat{f}_i - f^*_i |^2 = \frac{1}{H^2}\left(1 - \frac{H^2}{\nu}\right).
		\]
\end{fact}

\begin{remark}
	The condition for the Darboux transform of a cmc surfaces to be again a cmc surface can also be given in terms of linear conserved quantities \cite[Theorem~4.3]{burstall_discrete_2014}.
\end{remark}

\subsection{Discrete polarised curves and the monodromy problem}
We now briefly review the integrable theory of discrete polarised curves (see \cite{cho_periodic_nodate} for a detailed explanation).
For some simply-connected discrete domain $\Sigma^1 \subset \mathbb{Z}$, let $\mu$ be a nowhere vanishing function on the unoriented edges of $\Sigma^1$ so that $(\Sigma^1, \frac{1}{\mu})$ is a polarised domain.
Discrete maps $f : (\Sigma^1, \frac{1}{\mu}) \to \mathbb{H}$, or their lifts $F: (\Sigma^1, \frac{1}{\mu}) \to \mathbb{HP}^1$, are called discrete polarised curves.
Then the \emph{dual curve} $f^* : (\Sigma^1, \frac{1}{\mu}) \to \mathbb{H}$ of $f$ can be defined via
	\[
		\dif{f}^*_{ij} = \frac{1}{\mu_{ij}}\dif{f}_{ij}^{-1}.
	\]

The integrable structure of discrete polarised curves can again be described in terms of families of (flat) connections: putting $D(\lambda)_{ji} : \{i\} \times \mathbb{H}^2 \to  \{j\} \times \mathbb{H}^2$ with $\lambda \in \mathbb{R}$ defined on every edge $(ij)$ via
	\begin{equation}\label{eqn:flat1curve}
		D(\lambda)_{ji} := \Id_{ji} + \begin{pmatrix} 0 & \dif{f}_{ij} \\ \lambda \dif{f}^*_{ij} & 0 \end{pmatrix},
	\end{equation}
let $r(\lambda)_{ji} : \{i\} \times \mathbb{H}^2 \to  \{j\} \times \mathbb{H}^2$ be defined as
	\begin{equation}\label{eqn:flat2curve}
		r(\lambda)_{ji} := (G \bullet D(\lambda))_{ji} = G_j D(\lambda)_{ji} G_i
	\end{equation}
where $G = (e \,\, \psi)$ so that
	\begin{equation} \label{eqn:flatRepcurve}
		\begin{aligned}
			r(\lambda)_{ji} &= \Id_{ji} + \lambda \eta_{ji}
				:= \Id_{ji} + \lambda
					\begin{pmatrix}
						f_j \dif{f}^*_{ij} & -f_j \dif{f}^*_{ij} f_i \\
						\dif{f}^*_{ij} & -\dif{f}^*_{ij} f_i
					\end{pmatrix}\\
				&= \pi_i + \frac{\mu_{ij}-\lambda}{\mu_{ij}} \pi_j.
		\end{aligned}
	\end{equation}
    
The projective transformation $r(\lambda)^P_{ji} : \{i\} \times \mathbb{HP}^1 \to  \{j\} \times \mathbb{HP}^1$ induced by $r(\lambda)_{ji}$ defines a connection on the trivial bundle $\underline{\mathbb{HP}}^1$, called the \emph{(flat) connection associated to $f$}.

Darboux transformations of discrete polarised curves are given in terms of a cross-ratios condition \cite[Definition~3.3]{cho_periodic_nodate} (see also \cite[Definition~3.1]{cho_discrete_2021-1}):
\begin{definition}\label{def:curveDarb}
	Two discrete polarised curves $f, \hat{f}: (\Sigma^1, \frac{1}{\mu}) \to \mathbb{H}$ are called a \emph{Darboux pair} with spectral parameter $\nu$ if
		\[
			\cratio(f, f_j, \hat{f}_j, \hat{f}_i) = \frac{\nu}{\mu_{ij}}
		\]
	on every edge $(ij)$.
	In such case, one curve is called a \emph{Darboux transform} of the other curve with spectral parameter $\mu$.
\end{definition}
The Darboux transformations of discrete polarised curves can be given a characterisation in terms of the associated (flat) connections:
\begin{fact}
	Let $F: (\Sigma^1, \frac{1}{\mu}) \to \mathbb{HP}^1$ be a discrete polarised curve with associated connection $r(\lambda)^P$.
	Then $\hat{F}: (\Sigma^1, \frac{1}{\mu}) \to \mathbb{HP}^1$ is a Darboux transform of $L$ with spectral parameter $\nu$ if and only if $\hat{F}$ is $r(\nu)^P$--parallel.
\end{fact}

For explicit calculations, the $r(\nu)^P$--parallel condition \eqref{eqn:parallel} for Darboux transformations can be reformulated as follows.
Suppose that $\hat{F}$ is $r(\nu)^P$--parallel so that for some $\hat{\psi} \in \Gamma \hat{F}$, we have
	\[
		r(\nu)_{ji} \hat{\psi}_i = \hat{\psi}_j,
	\]
on any edge $(ij)$.
Now, if we define $\phi := G^{-1} \hat{\psi}$, then we have
	\[
		D(\nu)_{ji} \phi_i = G^{-1}_j r(\nu)_{ji} \hat{\psi}_i = \phi_j.
	\]
Therefore, if $\phi =: \begin{psmallmatrix} a \\ b \end{psmallmatrix}$ for some $a, b : \Sigma^1 \to \mathbb{H}$, then $a$ and $b$ must satisfy
	\begin{equation}\label{eqn:riccati2}
		\begin{pmatrix} \dif{a}_{ij} \\ \dif{b}_{ij} \end{pmatrix} = - \begin{pmatrix} \dif{f}_{ij} b_i \\ \nu \dif{f}^*_{ij} a_i \end{pmatrix},
	\end{equation}
and the new Darboux transform $\hat{f}$ is then given by $\hat{f} = f + a b^{-1}$.
	
To review the monodromy problem for Darboux transformations of discrete polarised curves \cite[\S~3.3]{cho_periodic_nodate}, suppose that $F, \hat{F}: (\Sigma^1, \frac{1}{\mu}) \to \mathbb{HP}^1$ are a Darboux pair of curves, and further assume that $F$ is an $M$-periodic curve, that is,
	\[
		F_m = F_{m + M}
	\]
for all $m \in \Sigma^1$.
If the \emph{monodromy matrix} $\mathcal{M}_{r, \nu}$ is defined by
	\[
		\mathcal{M}_{r, \nu} := \prod_{\iota = m}^{m+M-1} r(\nu)_{(\iota, \iota+1)},
	\]
then the eigenvectors $\phi$ of $\mathcal{M}_{r, \nu}$ correspond to \emph{parallel sections with multipliers}, that is, for $\hat{F} := \Span_\mathbb{H}\{ \phi \}$, we have
	\[
		\hat{F}_m = \hat{F}_{m + M}
	\]
where $\hat{F}$ is $r(\nu)^P$--parallel.
A spectral parameter $\nu$ is called a \emph{resonance point} if every parallel section is a parallel section with multiplier.
In other words, using a resonance point as a spectral parameter always gives closed Darboux transforms, regardless of the choice of the initial condition.

\section{Monodromy of discrete isothermic surfaces}\label{sec:monodromy}
The integrable structure of discrete isothermic surfaces and that of discrete polarised curves closely resemble each other.
In this section, we formalise the resemblance and interpret the integrability of discrete isothermic surfaces in terms of transformations of discrete polarised curves.
The reduction to curves will be crucial to examining the monodromy of Darboux transformations of discrete isothermic surfaces.
We note here that such reduction for the semi-discrete case has been considered in \cite{burstall_semi-discrete_2016}, where semi-discrete isothermic surfaces \cite{muller_semi-discrete_2013} are characterised as successive Darboux transforms of smooth polarised curves.

\subsection{Discrete isothermic surfaces in terms of discrete polarised curves}
Since we are interested in discrete curvature lines, let $\Sigma^2 \subset \mathbb{Z}^2$ be a discrete rectangular domain, that is, $\Sigma^2 = \Sigma^1 \times \tilde\Sigma^1$ for simply-connected domains $\Sigma^1, \tilde{\Sigma}^1 \subset \mathbb{Z}$.

Suppose that $f : \Sigma^2 \to \mathbb{H}$ is a discrete isothermic surface with cross-ratios factorising function $\mu$.
Viewing $\mu$ as a function defined on the unoriented edges of $\Sigma^2$, we can consider the pair $(\Sigma^2, \frac{1}{\mu})$ as a \emph{discrete polarised domain}, and view discrete isothermic surfaces as defined on such polarised domains.
Therefore, given a discrete isothermic surface on a polarised domain, all of its Christoffel transforms and Darboux transforms are defined on the same polarised domain.

Denote by $f^0 : \Sigma^1 \to \mathbb{H}$ any discrete curvature line of $f$ defined on the corresponding simply-connected domain $\Sigma^1$.
Since the cross-ratios factorising function $\mu$ serves as the polarisation of the domain $\Sigma^2$, we can consider the restriction of $\mu$ to the unoriented edges of $\Sigma^1$; therefore, we see that $f^0 : (\Sigma^1, \frac{1}{\mu}) \to \mathbb{H}$ is a discrete polarised curve.
In particular, we have that $\mu$ is an edge-labeling; thus every discrete curvature line in a family is defined on the same polarised domain $(\Sigma^1, \frac{1}{\mu})$.

Take any two neighbouring discrete curvature lines $f^0, f^1$ of $f$: without loss of generality, let $f^0 = f_{m, n_0}$ and $f^1 = f_{m, n_0+1}$ for some fixed $n_0$.
Then since $\mu$ is an edge-labeling, we can denote the common polarisation on the edges $((m,n_0), (m, n_0 + 1))$ for any $m$ as $\mu_{(n_0, n_0+1)}$.
The definition of discrete isothermic surfaces in Definition~\ref{def:discIsothermic} tells us that on any edge $(m, m+1)$, we have
	\[
		\cratio(f_{m, n_0}, f_{m + 1, n_0}, f_{m + 1, n_0 + 1}, f_{m, n_0 + 1}) = \frac{\mu_{(n_0, n_0+1)}}{\mu_{(m, m+1)}}.
	\]

Thus, the definition of Darboux transformations of discrete polarised curves in Definition~\ref{def:curveDarb} allows us to deduce the following:
\begin{proposition}
	Any two neighbouring discrete curvature lines of a discrete isothermic surface form a Darboux pair of discrete polarised curves.
\end{proposition}
\begin{remark}\label{rem:freedom}
	Note that the cross-ratios factorising functions are well-defined only up to constant multiplication.
	Thus, when restricting the polarised domain of a discrete isothermic surface to the polarised domain of a curvature line, there is a freedom of choice in the scaling of the polarisation.
\end{remark}

In a similar manner, for a Darboux pair of discrete isothermic surfaces $f, \hat{f} : (\Sigma^2, \frac{1}{\mu}) \to \mathbb{H}$, consider a corresponding pair of discrete curvature lines $f^0, \hat{f}^0 : (\Sigma^1, \frac{1}{\mu}) \to \mathbb{H}$ defined on a restricted discrete polarised domain.
Then we have for any edge $(ij) \subset \Sigma^1$,
	\[
		\cratio(f^0_i, f^0_j, \hat{f}^0_j, \hat{f}^0_i) = \frac{\nu}{\mu_{ij}}.
	\]
Thus the definition of Darboux transformations of discrete isothermic surfaces in Definition~\ref{def:darb} directly yield the following observation:
\begin{proposition}
	Any two corresponding discrete curvature lines of a discrete Darboux pair of isothermic surfaces form a Darboux pair of discrete polarised curves.
\end{proposition}

Furthermore, the well-definedness of Darboux transformations of discrete isothermic surfaces implies the permutability of Darboux transformations of discrete polarised curves:
\begin{theorem}\label{thm:perm}
	Let $f : (\Sigma^1, \frac{1}{\mu}) \to \mathbb{H}$ be a polarised curve, with Darboux transforms $f_1, f_2 : (\Sigma^1, \frac{1}{\mu}) \to \mathbb{H}$ of $f$ with respective spectral parameter $\nu_1, \nu_2$.
	Then there exists a fourth discrete polarised curve $f_{12} : (\Sigma^1, \frac{1}{\mu}) \to \mathbb{H}$ that is a simultaneous Darboux transform of $f_1$ and $f_2$ with respective spectral parameter $\nu_2$ and $\nu_1$.
\end{theorem}

\subsection{Reduction of the monodromy problem}
So far, we have obtained the reduction of the integrable structure of discrete isothermic surfaces and its Darboux transformations into a lattice of Darboux transformations of discrete polarised curves.
Such change in perspective is crucial for our consideration of the monodromy problem of Darboux transforms of discrete isothermic surfaces.

Thus, let $f$ be a discrete isothermic surface with a period $M$; without loss of generality, we will assume that $f$ is periodic in the $m$-curvature line direction so that $f_{m,n} = f_{m + M, n}$ for some $M \in \mathbb{N}$.

Now suppose that $\hat{f}$ is a Darboux transformation of $f$ with spectral parameter $\nu$.
Let us further assume that one of the $m$-curvature lines of $\hat{f}$ is periodic with the same period as $f$, i.e.\ 
	\[
		\hat{f}_{m,n_0} = \hat{f}_{m + M,n_0}
	\]
for some fixed $n_0$.
Then from the definition of Darboux transformations, the cross-ratios condition
	\[
		\cratio(f_{m,n_0}, f_{m,n_0 + 1}, \hat{f}_{m,n_0 + 1}, \hat{f}_{m,n_0}) = \frac{\nu}{\mu_{(n_0, n_0+1)}}
	\]
for any $m$ implies that $\hat{f}_{m, n_0+1}$ must be determined uniquely from the values of
	\[
		f_{m, n_0}, \quad f_{m, n_0 +1}, \quad  \hat{f}_{m,n_0}, \quad \nu, \quad \mu_{(n_0, n_0+1)}.
	\]
However, all of these values are $M$-periodic in the $m$-direction; thus, we have that $\hat{f}_{m,n_0+1}$ must also be periodic, that is,
	\[
		\hat{f}_{m, n_0+1} = \hat{f}_{m + M, n_0 +1}.
	\]
Propogating to the entire domain, we conclude that the entire surface is periodic, i.e., for any $n$,
	\[
		\hat{f}_{m, n} = \hat{f}_{m + M, n}.
	\]
We have thus proved that the monodromy problem of Darboux transformations of discrete isothermic surfaces can be reduced to a monodromy problem of Darboux transformations of a single discrete curvature line:
\begin{theorem}\label{thm:reduction}
	Let $f : (\Sigma^2, \frac{1}{\mu}) \to \mathbb{H}$ be a discrete isothermic surface with one family of periodic curvature lines.
	If $\hat{f} : (\Sigma^2, \frac{1}{\mu}) \to \mathbb{H}$ is a Darboux transform of $f$, then the corresponding family of curvature lines of $\hat{f}$ are also periodic with the same period if and only if one of the corresponding family of curvature lines is periodic.
\end{theorem}

Therefore, if $f$ is a discrete isothermic surfaces with periodic $m$-curvature lines, we can obtain new periodic isothermic surfaces as follows:
\begin{enumerate}
	\item Choose any base point $f_{m_0, n_0}$ of $f$ and the $m$-curvature line $f_{m, n_0}$ through that point.
	\item Viewing $f_{m, n_0}$ as a closed discrete polarised curve, find the spectral parameter $\nu$ and the initial condition $\hat{f}_{m_0, n_0}$ to obtain closed Darboux transform $\hat{f}_{m, n_0}$ (as a polarised curve).
	\item Use the spectral parameter $\nu$ and the initial condition $\hat{f}_{m_0, n_0}$ to obtain a new periodic discrete isothermic surface $\hat{f}$.
\end{enumerate}
In particular, if we have that all curvature lines of a discrete isothermic surface are periodic so that $f$ is a discrete torus, then we can also obtain new discrete isothermic tori as follows:
\begin{enumerate}
	\item Choose any base point $f_{m_0, n_0}$ of $f$, and take the $m$-curvature line $f_{m, n_0}$ and $n$-curvature line $f_{m_0, n}$ through that point.
	\item Viewing both curves as closed discrete polarised curves, find a spectral parameter $\nu$ and an initial condition $\hat{f}_{m_0, n_0}$ to obtain simultaneously closed Darboux transforms $\hat{f}_{m, n_0}$ and $\hat{f}_{m_0, n}$ of $f_{m, n_0}$ and $f_{m_0, n}$, respectively.
	\item Use the spectral parameter $\nu$ and the initial condition $\hat{f}_{m_0, n_0}$ to obtain a new periodic discrete isothermic torus $\hat{f}$.
\end{enumerate}

The main obstacle for the construction of new discrete isothermic tori via Darboux transformations is the second step.

\section{Periodic discrete isothermic surfaces}\label{sec:examples}
In this section, we put our theory to test to obtain new examples of periodic discrete isothermic surfaces.
In particular, we will consider the Darboux transforms of surfaces of revolution in $\mathbb{R}^3$, and obtain discrete analogues of isothermic bubbletons (see, for example, \cite[\S~2.2]{cho_generalised_2022}) via Darboux transformations.
By restricting to the case of discrete circular cylinders, we will also obtain the discrete analogues of cmc cylinders called bubbletons \cite{sterling_existence_1993}.
Finally, considering the Darboux transforms of homogeneous tori in $\mathbb{S}^3$, we will construct the discrete analogues of isothermic tori found in \cite{bernstein_non-special_2001} (see also \cite[Example~5.4.25]{hertrich-jeromin_introduction_2003}).

\subsection{Closed Darboux transforms of discrete circles -- revisited}\label{sect:circle}
The commonground for all the cases we consider in the section is the fact that the closed curvature lines will be multiply-covered discrete circles with constant polarisation.
Thus we first review the closure conditions for the Darboux transforms of multiply-covered discrete circles, where we generalise the calculations from \cite[Example~3.8]{cho_periodic_nodate} to allow for any radius and any constant polarisation.
We note here that we will only consider the discrete circles constructed via uniform sampling of the smooth circle, to obtain explicit results.

Let $f : (\Sigma^1, \frac{1}{\mu}) \to \mathbb{R}^2 \subset \mathbb{R}^4 \cong \Im \mathbb{H}$ be a $\rho$--fold cover of the discrete circle given by
	\[
		f_m := \jj r e^{\frac{2 \pi \ii}{M} m}
	\]
for some $r \in \mathbb{R}_+$ and $M \in \mathbb{N}$.
Throughout the paper, we will assume that $M > 2$.
Give $f$ a constant polarisation which we write as
	\[
		\frac{1}{\mu_{ij}} = \alpha |1 - e^\frac{2\pi \ii}{M}|^2
	\]
for some non-zero real constant $\alpha$.
Note that for $\alpha = r^2$, then we have
	\[
		\frac{1}{\mu_{ij}} = r^2 |1 - e^\frac{2\pi \ii}{M}|^2 = |f_i - f_j|^2 = |{\dif{f}_{ij}}|^2
	\]
so that $f$ is polarised by arc-length.
Then we obtain that $f$ is an $\rho M$--periodic discrete polarised curve for any $\rho \in \mathbb{N}$.

To consider the Darboux transforms, let $\hat{F} = \hat{\psi} \mathbb{H}$ be $r(\nu)^P$--parallel, and $\phi = \begin{psmallmatrix} a \\ b \end{psmallmatrix}$ where $\phi = G^{-1} \hat{\psi}$.
Rewriting \eqref{eqn:riccati2}, we have that on any three consecutive vertices $(ijp)$, $a$ must satisfy the recurrence relation
	\[
		\nu \dif{f}_{ij}^* a_i - (\dif{f}_{ij})^{-1} \dif{a}_{ij} + (\dif{f}_{jp})^{-1} \dif{a}_{jp} = 0.
	\]
In the case of discrete circles under consideration, the recurrence relation reads
	\[
		e^\frac{2\pi \ii}{M} a_p - (1 + e^\frac{2\pi \ii}{M}) a_j + (1 - \nu \alpha |1 - e^\frac{2\pi \ii}{M}|^2) a_i = 0.
	\]
If we set $a =: a_0 + \jj a_1$, for some functions $a_0, a_1 : \Sigma^2 \to \mathbb{C}$, we can solve the recurrence equation and obtain
	\begin{equation}\label{eqn:explicit}
		a = a_0^+c_0  + a_0^-c_1  + \jj(a_1^+c_2  + a_1^-c_3 )
	\end{equation}
for some constants $c_0, c_1, c_2, c_3 \in \mathbb{C}$, and
	\begin{equation}\begin{aligned}\label{eqn:explicit2}
		a_{0,m}^\pm &= \left(\tfrac{1}{2}(e^{-\frac{2\pi \ii}{M}}(1 \pm s) + (1 \mp s))\right)^m \\
		a_{1,m}^\pm &= \left(\tfrac{1}{2}(e^{\frac{2\pi \ii}{M}}(1 \mp s) + (1 \pm s))\right)^m
	\end{aligned}\end{equation}
where $s := \sqrt{1 - 4 \nu \alpha}$.
	
To consider the monodromy, we notice that $\phi = \begin{psmallmatrix}a \\ b \end{psmallmatrix}$ is an eigenvector of the monodromy matrix if and only if
	\begin{equation}\label{eqn:phipar}
		\phi_m h = \phi_{m + \rho M}
	\end{equation}
for some $h \in \mathbb{H}$.
However, due to the recurrence relation for $a,b$ in \eqref{eqn:riccati2}, we have that \eqref{eqn:phipar} holds if and only if
	\[
		a_m h = a_{m + \rho M}.
	\]

To investigate when the curve admits resonance points, we define $h_\iota^\pm$ so that
	\[
		a^\pm_{\iota,m} h_\iota^\pm = a^\pm_{\iota,m + \rho M}
	\]
for $\iota = 0, 1$, and look for conditions under which
	\begin{equation}\label{eqn:plzRP}
		h_0^+ = h_0^- = h_1^+ = h_1^-.
	\end{equation}
Since $a^\pm_{\iota,m + \rho M} = a^\pm_{\iota,m}a^\pm_{\iota,\rho M}$, we have $h_\iota^\pm = a^\pm_{\iota, \rho M}.$
However, note that
	\begin{align*}
		h_0^\pm = a^\pm_{0, \rho M} &= \left(\tfrac{1}{2}(e^{-\frac{2\pi \ii}{M}}(1 \pm s) + (1 \mp s))\right)^{\rho M} \\
			& = \left(\tfrac{1}{2}e^{-\frac{2\pi \ii}{M}} ((1 \pm s) + e^{\frac{2\pi \ii}{M}}(1 \mp s))\right)^{\rho M} \\
			&= a^\pm_{1, \rho M} = h_1^\pm.
	\end{align*}
Therefore, we have \eqref{eqn:plzRP} if and only if $h_1^+ = h_1^-$, that is,
	\[
		e^{2\pi \ii k} = \frac{h_1^+}{h_1^-} = \frac{a^+_{1,\rho M}}{a^-_{1,\rho M}} = \left(\frac{e^{\frac{2\pi \ii}{M}}(1 - s) + (1 + s)}{e^{\frac{2\pi \ii}{M}}(1 + s) + (1 - s)}\right)^{\rho M}
	\]
for some $k \in \mathbb{Z}$.
Thus, we have that $\nu$ is a resonance point of $f$ if and only if
	\begin{equation}\label{eqn:resonance}
		\nu = \frac{1}{4\alpha}\left(1 - \cot^2 \frac{\pi}{M} \tan^2 \frac{k \pi}{\rho M}\right).
	\end{equation}
We will assume without loss of generality that $\rho, k \in \mathbb{N}$.

We will now obtain the explicit parametrisations of closed Darboux transforms $\hat{f}$ of the discrete circle.
Taking a resonance point \eqref{eqn:resonance} as the spectral parameter, we have $s \in \mathbb{R}$ so that $a_{0}^\pm = \overline{a_{1}^\mp}$.
Thus, $a$ can be written as
	\begin{align*}
		a &= a_0^+c_0  + a_0^-c_1  + \overline{a_1^+}\jj c_2  + \overline{a_1^-} \jj c_3 = a_0^+ (c_0 + \jj c_3) + a_0^- (c_1 + \jj c_2) \\
		&=: a_0^+ c^+ + a_0^- c^-,
	\end{align*}
for some constants $c^\pm \in \mathbb{H}$.
We can further calculate that
	\[
		a^\pm_{0,m} = \left(e^{-\frac{\pi \ii}{M}} \cos \tfrac{\pi}{M} \sec \tfrac{k \pi}{\rho M}\right)^m e^{\mp\frac{k \pi \ii}{\rho M}m}
	\]
so that $a$ admits the following explicit form:
	\[
		a_m = \left(e^{-\frac{\pi \ii}{M}} \cos \tfrac{\pi}{M} \sec \tfrac{k \pi}{\rho M}\right)^m \left(e^{-\frac{k \pi \ii}{\rho M}m}c^+ + e^{\frac{k \pi \ii}{\rho M}m}c^-\right).
	\]

By \eqref{eqn:riccati2}, we have that
	\begin{align*}
		b_m &= -\dif{f}_{ij}^{-1} \dif{a}_{ij} \\ 
			&= \frac{\left(e^{-\frac{\pi \ii}{M}} \cos \tfrac{\pi}{M} \sec \tfrac{k \pi}{\rho M}\right)^m \jj \left( \sin{\frac{(\rho + k) \pi}{\rho M}}  e^{-\frac{k \pi \ii}{\rho M}m} c^+ + \sin{\frac{(\rho - k) \pi}{\rho M}} e^{\frac{k \pi \ii}{\rho M}m}  c^-\right)}{2r \sin{\tfrac{\pi}{M}} \cos{\tfrac{k \pi}{\rho M}}}.
	\end{align*}

Thus, every Darboux transform of a discrete circle $f$ is given by
	\begin{align*}
		\hat{f}_m &= f_m + a_m b_m^{-1}\\
			& = \jj r e^{\frac{2 \pi \ii}{M}m} + 2r \sin{\tfrac{\pi}{M}} \cos{\tfrac{k \pi}{\rho M}} e^{-\frac{\pi \ii}{M}m}\alpha_m \beta_m^{-1} e^{\frac{\pi \ii}{M}m}
	\end{align*}
where
	\begin{align*}
		\alpha_m &= e^{-\frac{k \pi \ii}{\rho M}m}c^+ + e^{\frac{k \pi \ii}{\rho M}m}c^-, \\
		\beta_m &= \jj \left(e^{-\frac{k \pi \ii}{\rho M}m} \sin{\tfrac{(\rho + k) \pi}{\rho M}} c^+ + e^{\frac{k \pi \ii}{\rho M}m} \sin{\tfrac{(\rho - k) \pi}{\rho M}} c^-\right)
	\end{align*}
for some constants $c^\pm \in \mathbb{H}$.

To find those Darboux transforms $\hat{f}$ living in the same $\jj \kk$--plane as $f$, we can assume without loss of generality that $c^+ = \jj$ and $c^- = \jj c_2$ for some $c_2 \in \mathbb{C}$.
Thus we obtain the explicit parametrisation of closed Darboux transforms of discrete circles living in the same $\jj \kk$--plane:
	\begin{equation}\label{eqn:yay}
			\hat{f}_m
				= - \jj r e^{\frac{2 \pi \ii}{M}m} \, \frac{e^{-\frac{k\pi \ii}{\rho M}m} \sin \frac{(\rho + k) \pi}{\rho M} c_2 + e^{\frac{k\pi \ii}{\rho M}m}   \sin \frac{(\rho - k) \pi}{\rho M}}{e^{-\frac{k\pi \ii}{\rho M}m} \sin \frac{(\rho - k) \pi}{\rho M} c_2 +  e^{\frac{k\pi \ii}{\rho M}m}\sin \frac{(\rho + k) \pi}{\rho M}}.
	\end{equation}

\subsection{Discrete isothermic cylinders}
Suppose we have a discrete surface of revolution in $\mathbb{R}^3 \cong \Im \mathbb{H}$ given by
	\[
		f_{m,n} := \ii q_n + \jj p_n e^{\frac{2\pi \ii}{M}m}
	\]
for some functions $p,q : \Sigma^1 \to \mathbb{R}$.
As calculated in \cite[\S~2]{burstall_discrete_2014}, we then have that for $i = (m,n)$,
	\begin{equation}\label{eqn:crsor}
		\cratio(f_i, f_j, f_k, f_\ell) = -\frac{p_n p_{n+1}}{(p_n - p_{n+1})^2 + (q_n - q_{n+1})^2} |1 - e^{\frac{2 \pi \ii}{M}}|^2 = \frac{\mu_{i\ell}}{\mu_{ij}}.
	\end{equation}
Clearly, $f$ has $M$-periodic $m$-curvature lines.

To consider the monodromy of its Darboux transforms over $\rho$--fold cover in the $m$-direction, we use Theorem~\ref{thm:reduction} and consider a single $m$-curvature line, say $f_{m,n_0}$ for some fixed $n_0$.
By Remark~\ref{rem:freedom}, we have some freedom to choose the correct polarisation $\mu$ for the discrete polarised curve $f_{m,n_0}$ from \eqref{eqn:crsor}.
Thus we choose
	\[
		\frac{1}{\mu_{ij}} = |1 - e^{\frac{2 \pi \ii}{M}}|^2
	\]
so that any two neighbouring $m$-curvature lines will be a Darboux pair with spectral parameter $\mu_{i\ell}$ where
	\[
		\frac{1}{\mu_{i\ell}} = -\frac{(p_n - p_{n+1})^2 + (q_n - q_{n+1})^2}{p_n p_{n+1}}.
	\]
If we assume without loss of generality that $q_{n_0} = 0$, then we have
	\[
		f_{m,n_0} = \jj p_{n_0} e^{\frac{2\pi \ii}{M}m}.
	\]
Thus we have that $f_{m,n_0}$ is a discrete polarised circle with constant polarisation.

Using the discussions in Section~\ref{sect:circle}, if we choose the spectral parameter $\nu$ as the resonance point \eqref{eqn:resonance} with $\alpha = 1$, then any Darboux transform $\hat{f}_{m,n_0}$ of the discrete polarised curve $f_{m,n_0}$ will be periodic.
Therefore, Theorem~\ref{thm:reduction} tells us that any Darboux transform $\hat{f}$ of $f$ with spectral parameter $\nu$ satisfying
	\[
		\nu = \frac{1}{4}\left(1 - \cot^2 \frac{\pi}{M} \tan^2 \frac{k \pi}{\rho M}\right)
	\]
gives a new discrete isothermic cylinder with period $\rho M$.

\subsection{Discrete isothermic bubbletons via discrete circular cylinder}\label{sect:circCyliso}
Now suppose that $f$ is a discrete circular cylinder, given by
	\[
		f_{m,n} := \ii q_n + \jj e^{\frac{2\pi \ii}{M}m}.
	\]
Then this is a special case of the discrete surface of revolution, so that
	\[
		\cratio(f_i, f_j, f_k, f_\ell) = -\frac{1}{(q_n - q_{n+1})^2} |1 - e^{\frac{2 \pi \ii}{M}}|^2 = \frac{\mu_{i\ell}}{\mu_{ij}}.
	\]
We again take
	\[
		\frac{1}{\mu_{ij}} = |1 - e^{\frac{2 \pi \ii}{M}}|^2, \quad \frac{1}{\mu_{i\ell}} = -(q_n - q_{n+1})^2.
	\]
Then any Darboux transform $\hat{f}$ of $f$ with spectral parameter 
	\begin{equation}\label{eqn:cylRes}
		\nu = \frac{1}{4}\left(1 - \cot^2 \frac{\pi}{M} \tan^2 \frac{k \pi}{\rho M}\right)
	\end{equation}
gives a discrete isothermic cylinder.

We now aim to obtain explicit parametrisations for such discrete cylinders: for this, we assume that
	\begin{enumerate}
		\item the spectral parameter $\nu$ is a resonance point \eqref{eqn:cylRes}, and
		\item $q_n = \frac{n}{N}$ for some $N \in \mathbb{Z}$.
	\end{enumerate}
We will solve for $a, b : \Sigma^2 \to \mathbb{H}$ so that $\hat{f} = f + a b^{-1}$, that is, $a$ and $b$ satisfies the linearised Riccati equation \eqref{eqn:riccati2} on any edge $(ij)$.
Since we have $f_{m,0} = \jj e^{\frac{2\pi \ii}{M}m}$, the discussions in Section~\ref{sect:circle} tells us that
	\begin{equation}\label{eqn:cylcirc}
		\begin{aligned}
			a_{m,0} &= \left(e^{-\frac{\pi \ii}{M}} \cos \tfrac{\pi}{M} \sec \tfrac{k \pi}{\rho M}\right)^m \alpha_m, \\
			b_{m,0} &= \frac{1}{2r}  \csc{\tfrac{\pi}{M}} \sec{\tfrac{k \pi}{\rho M}} \left(e^{-\frac{\pi \ii}{M}} \cos \tfrac{\pi}{M} \sec \tfrac{k \pi}{\rho M}\right)^m \beta_m
		\end{aligned}
	\end{equation}
where
	\begin{align*}
		\alpha_m &= e^{-\frac{k \pi \ii}{\rho M}m}c^+ + e^{\frac{k \pi \ii}{\rho M}m}c^-, \\
		\beta_m &= \jj \left(e^{-\frac{k \pi \ii}{\rho M}m} \sin{\tfrac{(\rho + k) \pi}{\rho M}} c^+ + e^{\frac{k \pi \ii}{\rho M}m} \sin{\tfrac{(\rho - k) \pi}{\rho M}}  c^-\right)
	\end{align*}
for some $c^\pm \in \mathbb{H}$.
%

Now, fixing any $m = m_0 \in \mathbb{Z}$, we will solve for $a, b$ over the line $(m_0, n)$.
Thus, we first need to solve for $a$ where $a$ satisfies the recurrence relation
	\[
		\nu \dif{f}_{i \ell}^* a_i - (\dif{f}_{i\ell})^{-1} \dif{a}_{i\ell} + (\dif{f}_{\ell q})^{-1} \dif{a}_{\ell q} = 0
	\]
over any three consecutive vertices $i = (m_0, n), \ell = (m_0, n+1), q = (m_0, n+2)$.
Under current assumptions, the recurrence relation reads
	\[
		a_q - 2 a_\ell + \left(1+ \frac{\nu}{N^2}\right) a_i = 0.
	\]
This allows us to solve for $a_{m_0,n}$ explicitly, and obtain that
	\[
		a_{m,n} = \left(1+ \tfrac{\sqrt{-\nu}}{N}\right)^n \gamma_m^+ + \left( 1- \tfrac{\sqrt{-\nu}}{N}\right)^n \gamma_m^-
	\]
where $\gamma^\pm$ are quaternionic-valued functions over $m$ alone.
Then the explicit expression for $b$ follows via \eqref{eqn:riccati2}:
	\[
		b_{m,n} = \ii \left(  \left(1+ \tfrac{\sqrt{-\nu}}{N}\right)^n \gamma_{m}^+ - \left(1 - \tfrac{\sqrt{-\nu}}{N}\right)^n \gamma_{m}^- \right) \sqrt{-\nu}.
	\]
To find $\gamma_m^\pm$, note that when $n = 0$, we have
	\[
		a_{m,0} = \gamma_{m}^+ +  \gamma_{m}^- \quad\text{and}\quad b_{m,0} = \ii (\gamma_{m}^+ -  \gamma_{m}^-)\sqrt{-\nu}
	\]
so that
	\[
		\gamma_m^\pm = \frac{1}{2}\left( a_{m,0} \mp \ii b_{m,0} \frac{1}{\sqrt{-\nu}}\right).
	\]
Using \eqref{eqn:cylcirc}, we find that
	\[
		\gamma_m^\pm = \frac{1}{2}\left(e^{-\frac{\pi \ii}{M}} \cos \tfrac{\pi}{M} \sec \tfrac{k \pi}{\rho M}\right)^m \left( \alpha_m \mp \frac{\ii}{2 \sqrt{-\nu}} \csc{\tfrac{\pi}{M}} \sec{\tfrac{k \pi}{\rho M}}  \beta_m \right).
	\]
Thus, we have found closed-form solutions for the functions $a$ and $b$, which in turn yields the explicit parametrisations of the discrete isothermic cylinder that are Darboux transforms of the discrete circular cylinder via
	\[
		\hat{f}_{m,n} = f_{m,n} + a_{m,n} b_{m,n}^{-1}.
	\]
	
\begin{figure}
	\centering
	\begin{minipage}{0.46\linewidth}
		\centering
		\includegraphics[width=\textwidth]{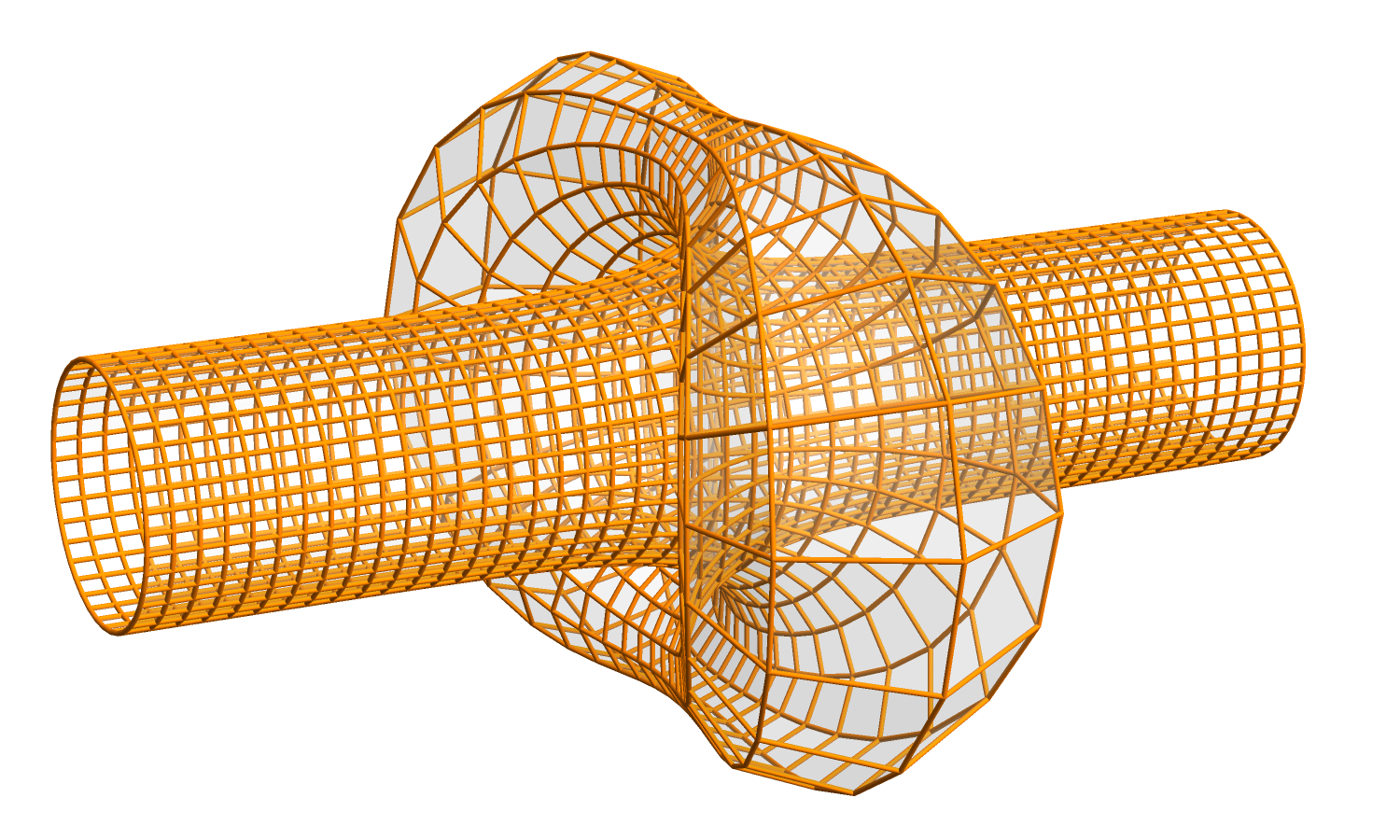}
	\end{minipage}
	\begin{minipage}{0.46\linewidth}
		\centering
		\includegraphics[width=\textwidth]{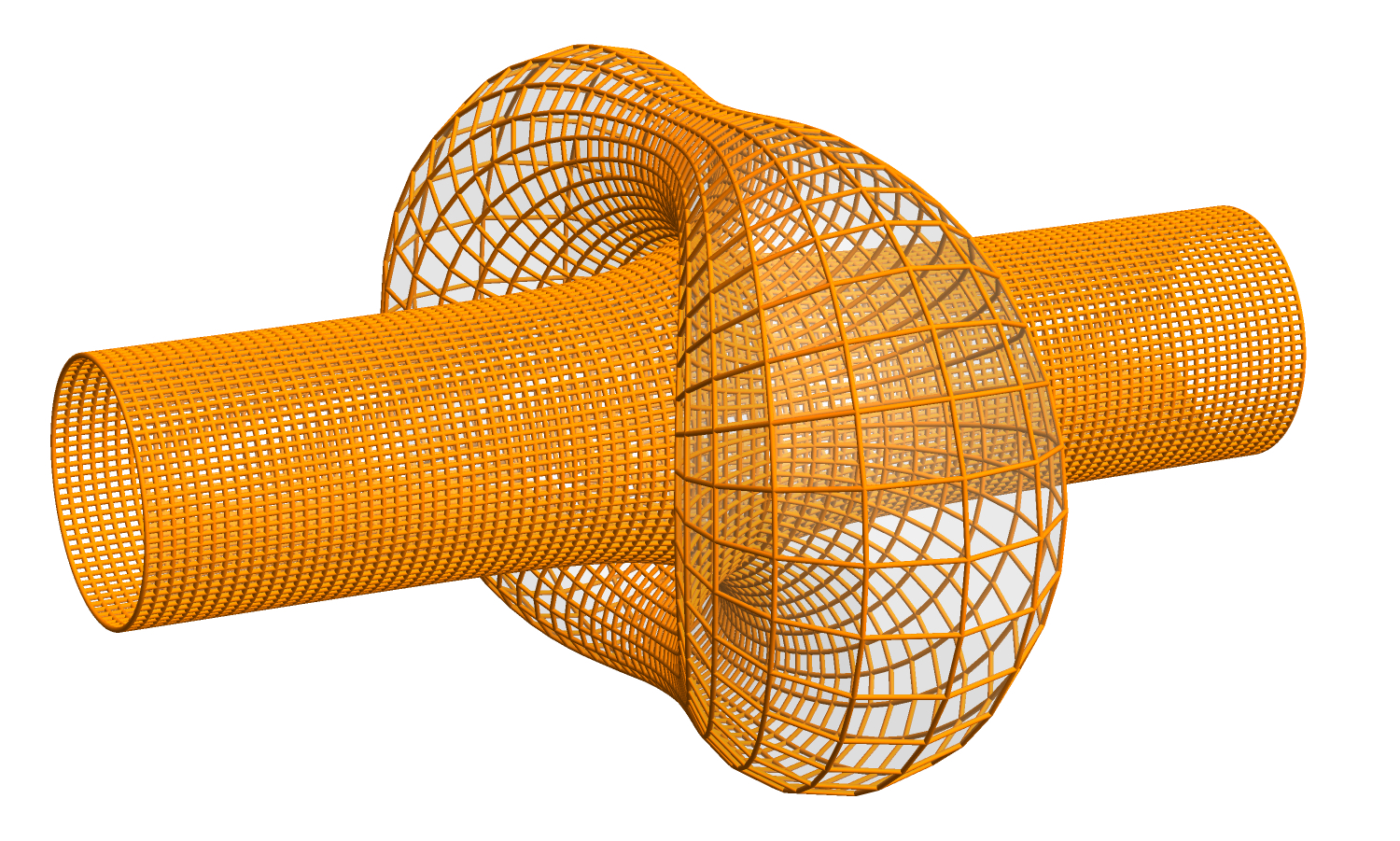}
	\end{minipage}
	\caption{Discrete isothermic bubbletons drawn using explicit parametrisation \eqref{eqn:expPar} with $k = 2$, $\rho = 1$, and $c_2 = -10$ with various number of subdivisions (on the left: $M = 40$, $N = 5$, on the right: $M = 80$, $N = 10$).}
	\label{fig:isoBubbleton}
\end{figure}

\begin{figure}
	\centering
	\begin{minipage}{0.46\linewidth}
		\centering
		\includegraphics[width=\textwidth]{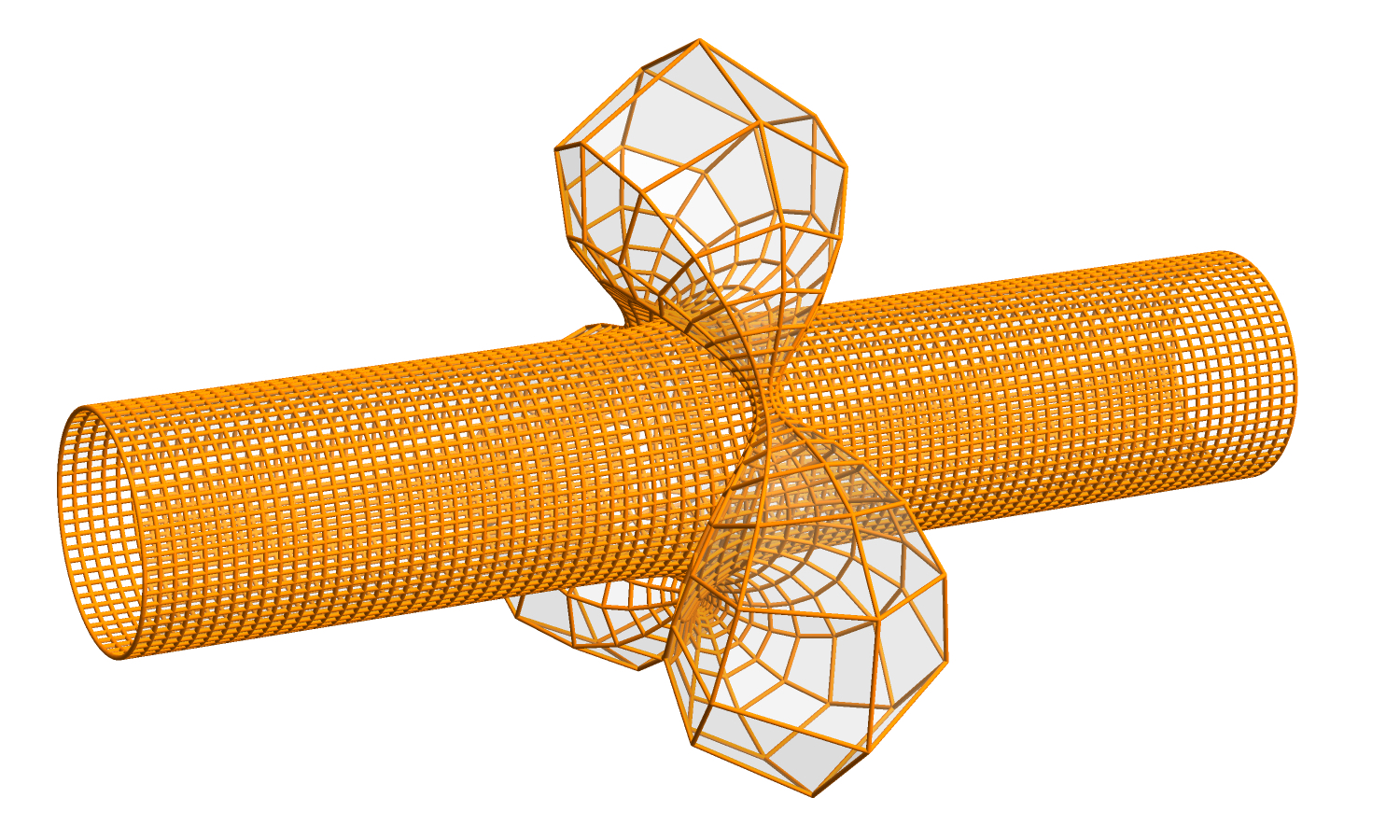}
	\end{minipage}
	\begin{minipage}{0.46\linewidth}
		\centering
		\includegraphics[width=\textwidth]{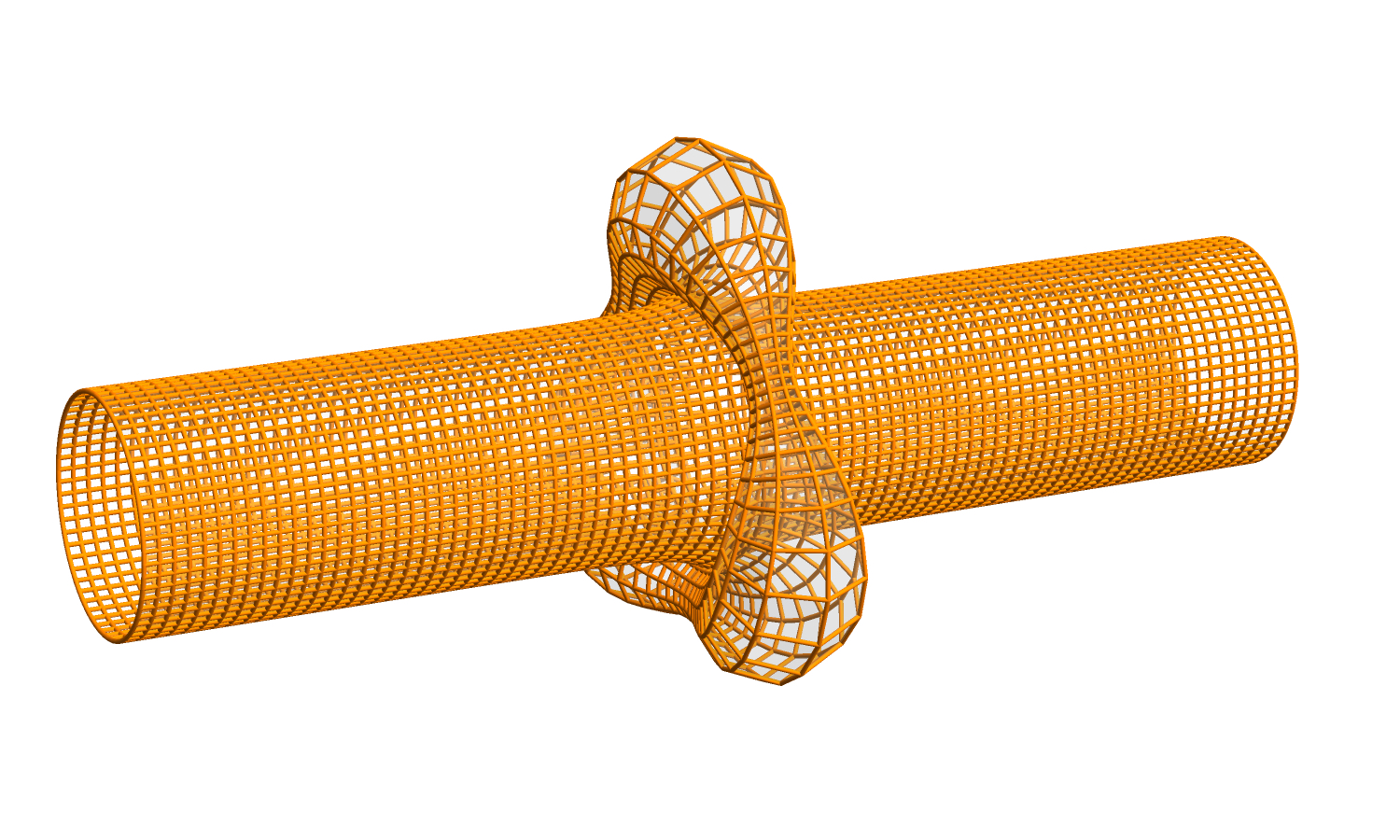}
	\end{minipage}
	\caption{Discrete isothermic bubbletons drawn using explicit parametrisation \eqref{eqn:expPar} with $k = 3$, $\rho = 1$, $M = 80$, and $N = 8$ with varying initial conditions (on the left: $c_2 = -4$, on the right: $c_2 = -8$).}
	\label{fig:isoBubbleton2}
\end{figure}

\begin{example}\label{exam:explicit}
For example, if we additionally assume that
	\begin{enumerate}
		\item the spectral parameter at the resonance point satisfies $\nu < 0$, and
		\item $\hat{f}_{m, 0}$ lives in the same $\jj \kk$--plane as $f_{m,0} = \jj e^{\frac{2\pi \ii}{M}m}$,
	\end{enumerate}
then without loss of generality we may set $c^+ = \jj$ and $c^- = \jj c_2$ for some $c_2 \in \mathbb{C}$.
Then we have
	\begin{align*}
		\alpha_m &= \jj \left(e^{\frac{k \pi \ii}{\rho M}m} + e^{-\frac{k \pi \ii}{\rho M}m}c_2 \right) =: \jj A_m \\ 
		\beta_m &= - \left(e^{\frac{k \pi \ii}{\rho M}m} \sin{\tfrac{(\rho + k) \pi}{\rho M}} + e^{-\frac{k \pi \ii}{\rho M}m} \sin{\tfrac{(\rho - k) \pi}{\rho M}}  c_2\right) =: -B_m,
	\end{align*}
so that $\gamma^\pm$ takes the following form:
	\[
		\gamma_m^\pm = \frac{1}{2}\left(e^{-\frac{\pi \ii}{M}} \cos \tfrac{\pi}{M} \sec \tfrac{k \pi}{\rho M}\right)^m \left( \jj A_m \pm \frac{\ii}{2 \sqrt{-\nu}} \csc{\tfrac{\pi}{M}} \sec{\tfrac{k \pi}{\rho M}}  B_m \right).
	\]
From this, we calculate that
	\begin{align*}
		|\gamma_m^\pm|^2 &= \frac{1}{4}\left(\cos \tfrac{\pi}{M} \sec \tfrac{k \pi}{\rho M}\right)^{2m} \left( |A_m|^2 - \tfrac{1}{4\nu} \csc^2{\tfrac{\pi}{M}} \sec^2{\tfrac{k \pi}{\rho M}}  |B_m|^2 \right),\\
		\gamma_m^+ \overline{\gamma_m^-} + \gamma_m^- \overline{\gamma_m^+}  &= \frac{1}{2}\left(\cos \tfrac{\pi}{M} \sec \tfrac{k \pi}{\rho M}\right)^{2m} \left( |A_m|^2 + \tfrac{1}{4\nu} \csc^2{\tfrac{\pi}{M}} \sec^2{\tfrac{k \pi}{\rho M}}  |B_m|^2 \right), \\
		\gamma_m^- \overline{\gamma_m^+} - \gamma_m^+ \overline{\gamma_m^-}  &= -\jj \frac{\ii}{2 \sqrt{-\nu}}\left(e^{\frac{\pi \ii}{M}} \cos \tfrac{\pi}{M} \sec \tfrac{k \pi}{\rho M}\right)^{2m} \csc{\tfrac{\pi}{M}} \sec{\tfrac{k \pi}{\rho M}} A_m \overline{B_m}.
	\end{align*}
Using these equations, we can write the following expressions explicitly:
	\begin{align*}
		|b_{m,n}|^2 &= 
			\begin{multlined}[t]
				-\nu \Big( \big(1+ \tfrac{\sqrt{-\nu}}{N}\big)^{2n} |\gamma_m^+|^2 + \big(1- \tfrac{\sqrt{-\nu}}{N}\big)^{2n} |\gamma_m^-|^2 \\
					\phantom{aaaaa} - \big(1+ \tfrac{\sqrt{-\nu}}{N}\big)^{n} \big(1- \tfrac{\sqrt{-\nu}}{N}\big)^{n} (\gamma_m^+ \overline{\gamma_m^-} + \gamma_m^- \overline{\gamma_m^+}) \Big),
			\end{multlined} \\
		a_{m,n} \overline{b_{m,n}} &=
			\begin{multlined}[t]
				-\sqrt{-\nu}  \Big( \big(1+ \tfrac{\sqrt{-\nu}}{N}\big)^{2n} |\gamma_m^+|^2 - \big(1- \tfrac{\sqrt{-\nu}}{N}\big)^{2n} |\gamma_m^-|^2 \\
					\phantom{aaaaa} + \big(1+ \tfrac{\sqrt{-\nu}}{N}\big)^{n} \big(1- \tfrac{\sqrt{-\nu}}{N}\big)^{n} (\gamma_m^- \overline{\gamma_m^+} - \gamma_m^+ \overline{\gamma_m^-}) \Big) \ii.
			\end{multlined}
	\end{align*}

Therefore, we obtain the closed-form parametrisation of the discrete isothermic cylinders $\hat{f}$ via
	\begin{equation}\label{eqn:expPar}
		\hat{f} = f + a b^{-1} = f + a \bar{b}|b|^{-2} =: f + \frac{1}{C}(T^0 + \jj T^1)
	\end{equation}
where
	\begin{align*}
		C_{m,n} &=
			\begin{multlined}[t]
				 2\big((N - \sqrt{-\nu})^n + (N + \sqrt{-\nu})^n\big)^2 |B_m|^2 \\
					+ \big((N - \sqrt{-\nu})^n - (N + \sqrt{-\nu})^n\big)^2 |A_m|^2(\cos \tfrac{2\pi}{M} - \cos \tfrac{2 k \pi}{\rho M})
			\end{multlined} \\
		T^0_{m,n} &=
				\frac{\ii}{\sqrt{-\nu}} \big((N - \sqrt{-\nu})^{2n} - (N + \sqrt{-\nu})^{2n}\big)
					\big(2|B_m|^2 + |A_m|^2(\cos \tfrac{2\pi}{M} - \cos \tfrac{2 k \pi}{\rho M})\big)\\
		T^1_{m,n} &=
				-16 e^{\frac{2 \pi \ii}{M}m} (N - \sqrt{-\nu})^n (N +\sqrt{-\nu})^n \sin \tfrac{\pi}{M} \cos \tfrac{k \pi}{\rho M} A_m \overline{B_m}.
	\end{align*}
See Figures~\ref{fig:isoBubbleton} and \ref{fig:isoBubbleton2} for examples.
\end{example}

\subsection{Discrete cmc bubbletons}\label{sect:circCylcmc}
Discrete circular cylinders are discrete cmc surfaces, which can be verified as follows: noting that
	\begin{align*}
		f^*_{m,n} - f^*_{m+1,n} &= \dif{f}^*_{ij} = \frac{1}{\mu_{ij}} \dif{f}_{ij}^{-1} = - \jj e^{\frac{2\pi \ii}{M}m} + \jj e^{\frac{2\pi \ii}{M}(m+1)} \\
		f^*_{m,n} - f^*_{m,n+1} &= \dif{f}^*_{i\ell} = \frac{1}{\mu_{i\ell}} \dif{f}_{i\ell}^{-1} = \ii q_n - \ii q_{n+1},
	\end{align*}
we take
	\[
		f^*_{0,0} = \ii q_0 -\jj
	\]
as the initial condition to obtain
	\[
		f^*_{m,n} = \ii q_n - \jj e^{\frac{2\pi \ii}{M}m}
	\]
as the dual surface.
Therefore, we have
	\[
		f^*_{m,n} - f_{m,n} = -2 \jj e^{\frac{2\pi \ii}{M}m},
	\]
so that
	\[
		|f^* - f|^2 = 4,
	\]
allowing us to conclude that $f$ is a discrete cmc surface with $H \equiv \frac{1}{2}$ by Definition~\ref{def:cmcDef}.

Thus, if we can choose the correct initial condition using Fact~\ref{fact:cmcCond}, then we will obtain new discrete cmc cylinders via Darboux transformations.
To find such initial condition, we first note that
	\[
		\cratio(f_i, f_j, f^*_j, f^*_i) = \frac{1}{4} |1 - e^{\frac{2\pi \ii}{M}}|^2 \quad\text{and}\quad \cratio(f_i, f_\ell, f^*_\ell, f^*_i) = -\frac{1}{4} (q_n - q_{n+1})^2
	\]
so that $f^*$ is a Darboux transform of $f$ with spectral parameter $H^2 = \frac{1}{4}$ (confirming Fact~\ref{fact:cmcDarb}).
Thus for some point $i \in \Sigma^2$, we must find $\hat{f}_i$ so that
	\[
		|\hat{f}_i - f^*_i|^2 = 4\left(1-\frac{1}{4\nu}\right),
	\]
which exists if and only if $\nu > \frac{1}{4}$ or $\nu < 0$.

\begin{figure}
	\centering
	\begin{minipage}{0.46\linewidth}
		\centering
		\includegraphics[width=0.9\textwidth]{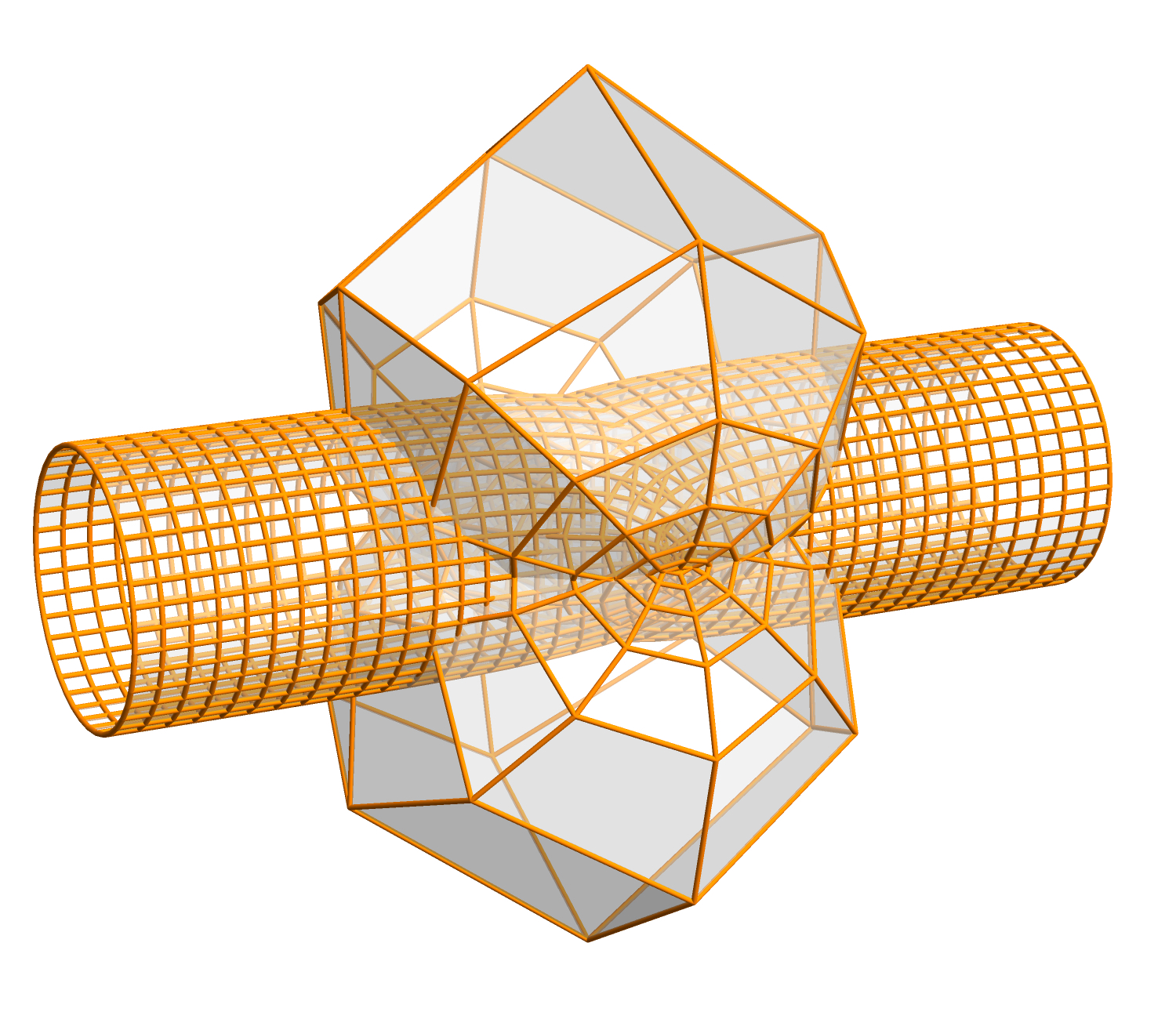}
	\end{minipage}
	\begin{minipage}{0.46\linewidth}
		\centering
		\includegraphics[width=0.9\textwidth]{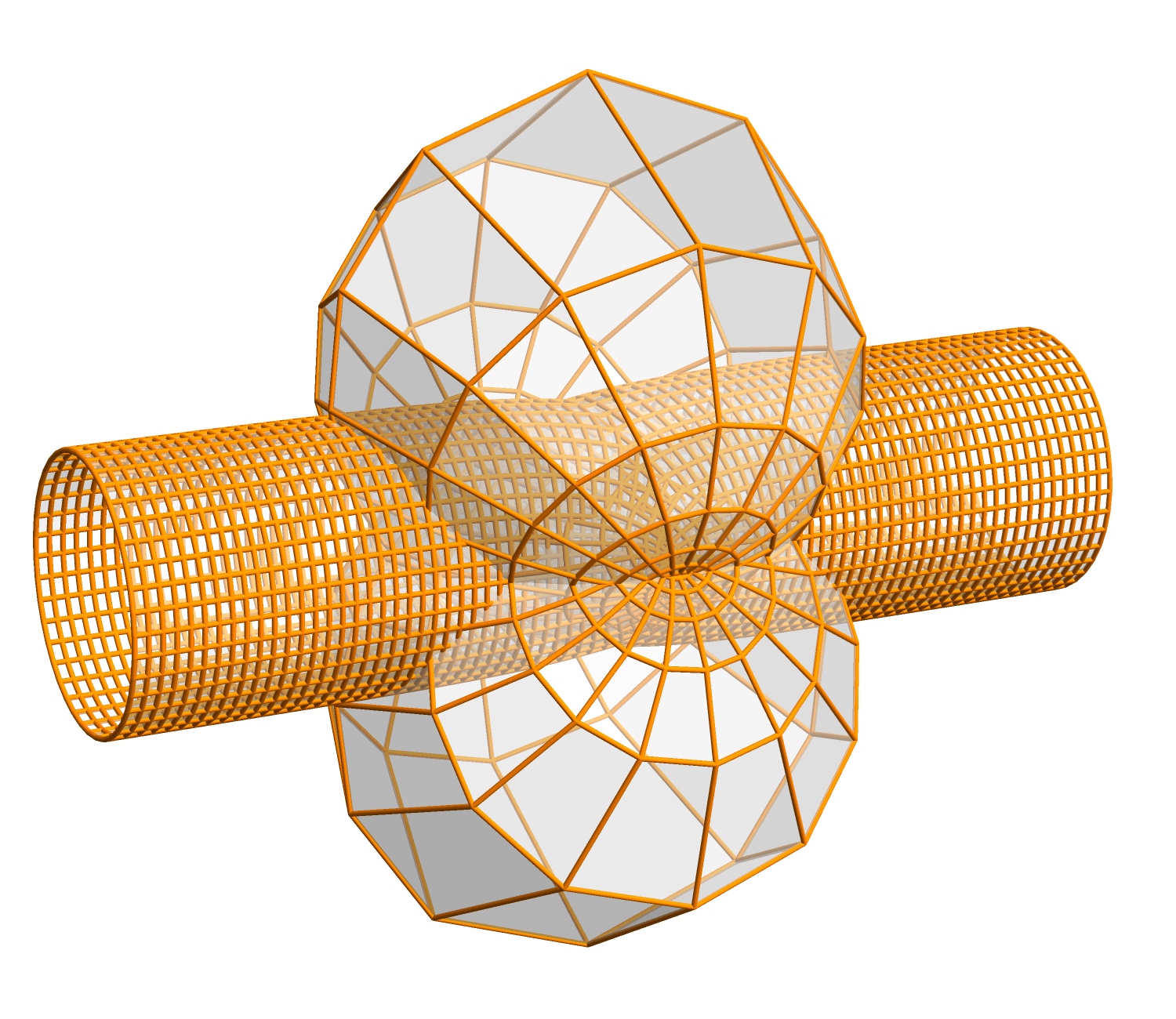}
	\end{minipage}
	\begin{minipage}{0.46\linewidth}
		\centering
		\includegraphics[width=0.9\textwidth]{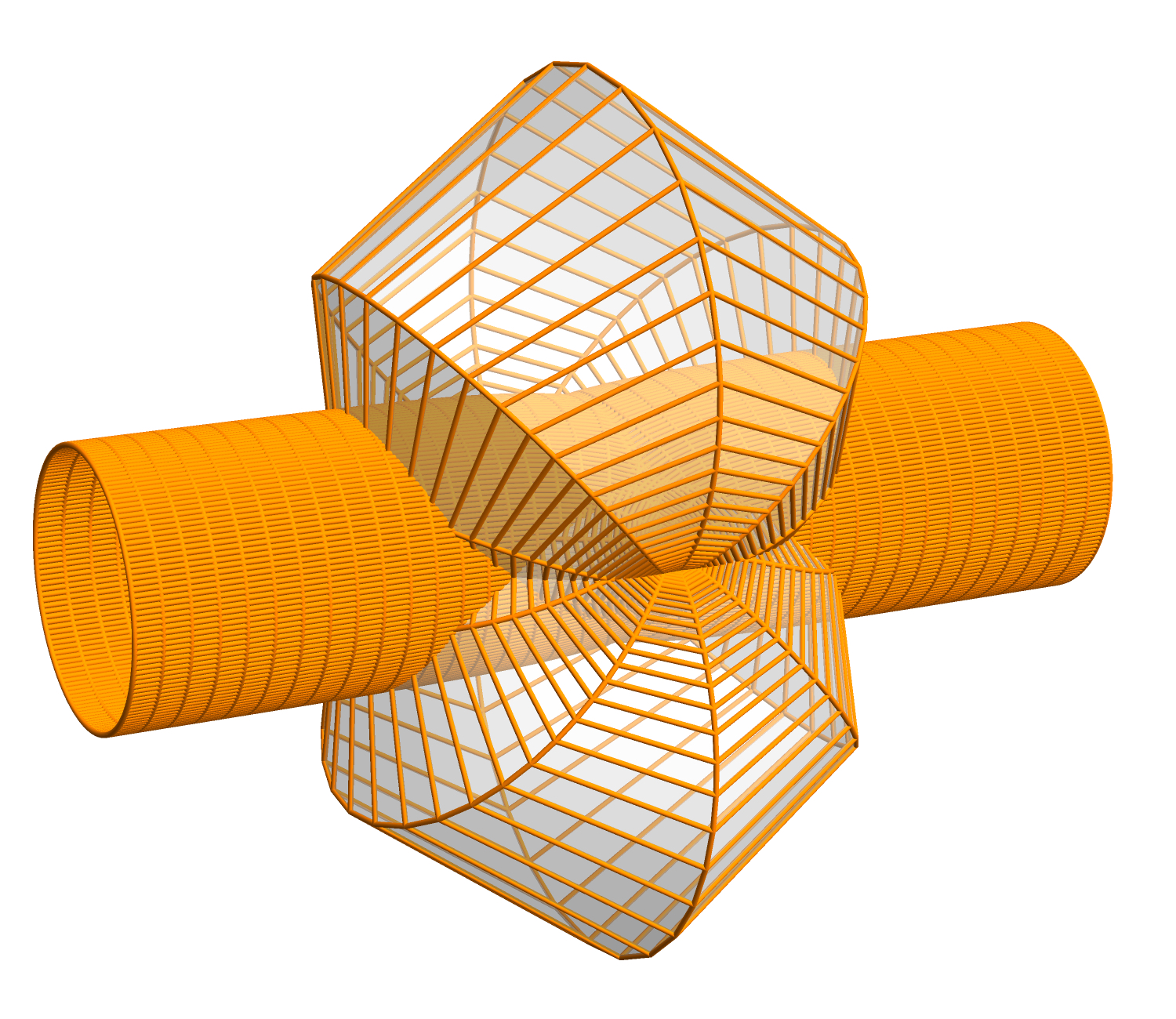}
	\end{minipage}
	\begin{minipage}{0.46\linewidth}
		\centering
		\includegraphics[width=0.9\textwidth]{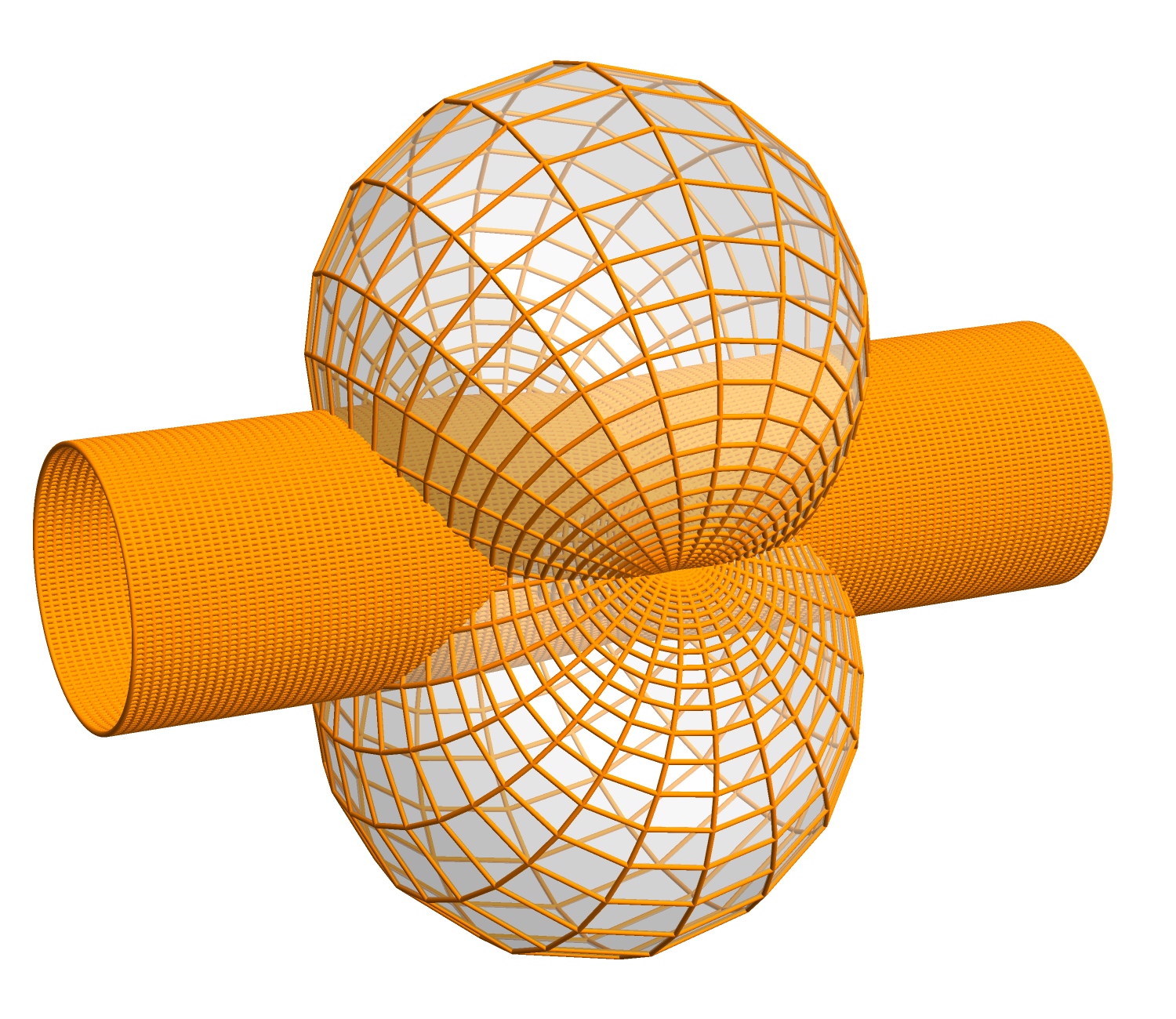}
	\end{minipage}
	\caption{Discrete cmc bubbletons drawn using explicit parametrisation \eqref{eqn:expPar} where $k = 2$, $\rho = 1$ with initial condition as in \eqref{eqn:cmcInit} at various number of subdivisions (on the top left: $M = 40$, $N = 5$; on the top right: $M = 40$, $N = 10$; on the bottom left: $M = 160$, $N = 5$; on the bottom right: $M = 160$, $N = 15$).}
	\label{fig:bubbleton}
\end{figure}

To investigate the conditions on $k$ and $\rho$ for the spectral parameter at a resonance point
	\[
		\nu = \frac{1}{4}\left(1 - \cot^2 \frac{\pi}{M} \tan^2 \frac{k \pi}{\rho M}\right)
	\]
to satisfy either $\nu > \frac{1}{4}$ or $\nu < 0$, first note that  $k = \rho$ implies $\nu = 0$; thus, we must have $k \neq \rho$.
Suppose that $0< k < \rho$.
On the inverval $u \in (0, \frac{\pi}{2})$, $\sin u$ is a strictly increasing positive function while $\cos u$ is a strictly decreasing positive function so that
	\[
		0 < \sin \frac{k \pi}{\rho M} < \sin \frac{\pi}{M} \quad\text{and}\quad 0 < \cos \frac{\pi}{M} < \cos \frac{k \pi}{\rho M}.
	\]
Thus, we deduce
	\[
		0 < \tan \frac{k \pi}{\rho M} < \tan \frac{\pi}{M},
	\]
that is,
	\[
		0 < \cot^2 \frac{\pi}{M} \tan^2 \frac{k \pi}{\rho M} < 1.
	\]
This implies that
	\[
		0 < \nu < \frac{1}{4}
	\]
so that we do not have discrete cmc cylinders appearing as Darboux transforms of the discrete circular cylinders when $k < \rho$.
This is the discrete analogue of the smooth case investigated in \cite{sterling_existence_1993} (see also \cite{cho_new_2022}).

On the other hand, if $k > \rho$, then for large enough $M$ so that $\frac{k}{\rho M} < \frac{1}{2}$, we have
	\[
		0 < \sin \frac{\pi}{M} < \sin \frac{k \pi}{\rho M} < 1 \quad\text{and}\quad 0 < \cos \frac{k \pi}{\rho M} < \cos \frac{\pi}{M}.
	\]
Thus, in this case
	\[
		0 < \tan \frac{\pi}{M} < \tan \frac{k \pi}{\rho M}
	\]
so that
	\[
		1 < \cot^2  \frac{\pi}{M} \tan^2 \frac{k \pi}{\rho M}.
	\]
In other words, we have $\nu < 0$.
Thus, when we choose the spectral parameter $\nu$ at the resonance point with $k > \rho$ with large enough $M$, we can find a correct initial condition to obtain another discrete cmc cylinder via Darboux transformation.
These are the discrete analogues of the bubbletons constructed in \cite{sterling_existence_1993}.

\begin{figure}
	\centering
	\begin{minipage}{0.46\linewidth}
		\centering
		\includegraphics[width=\textwidth]{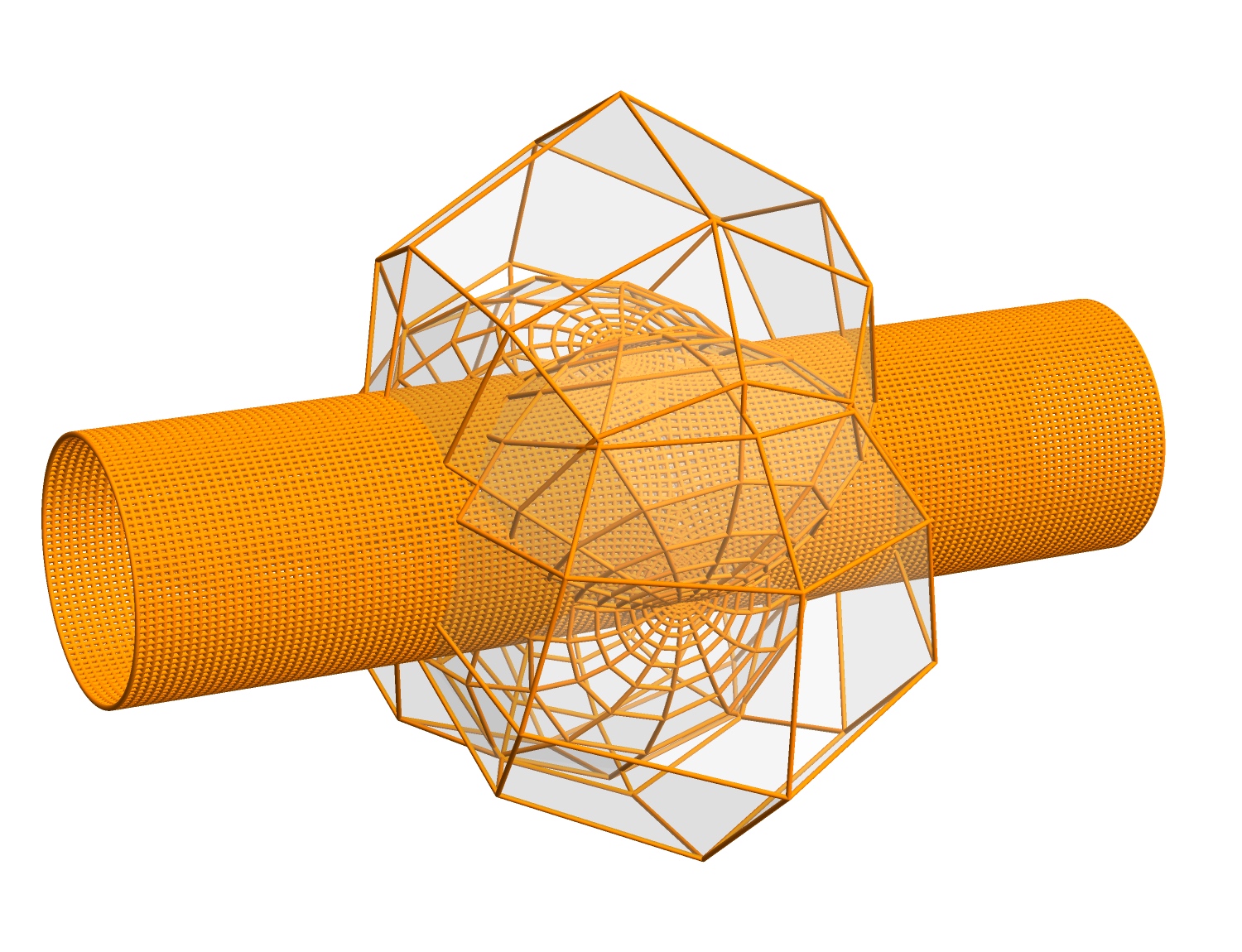}
	\end{minipage}
	\begin{minipage}{0.46\linewidth}
		\centering
		\includegraphics[width=\textwidth]{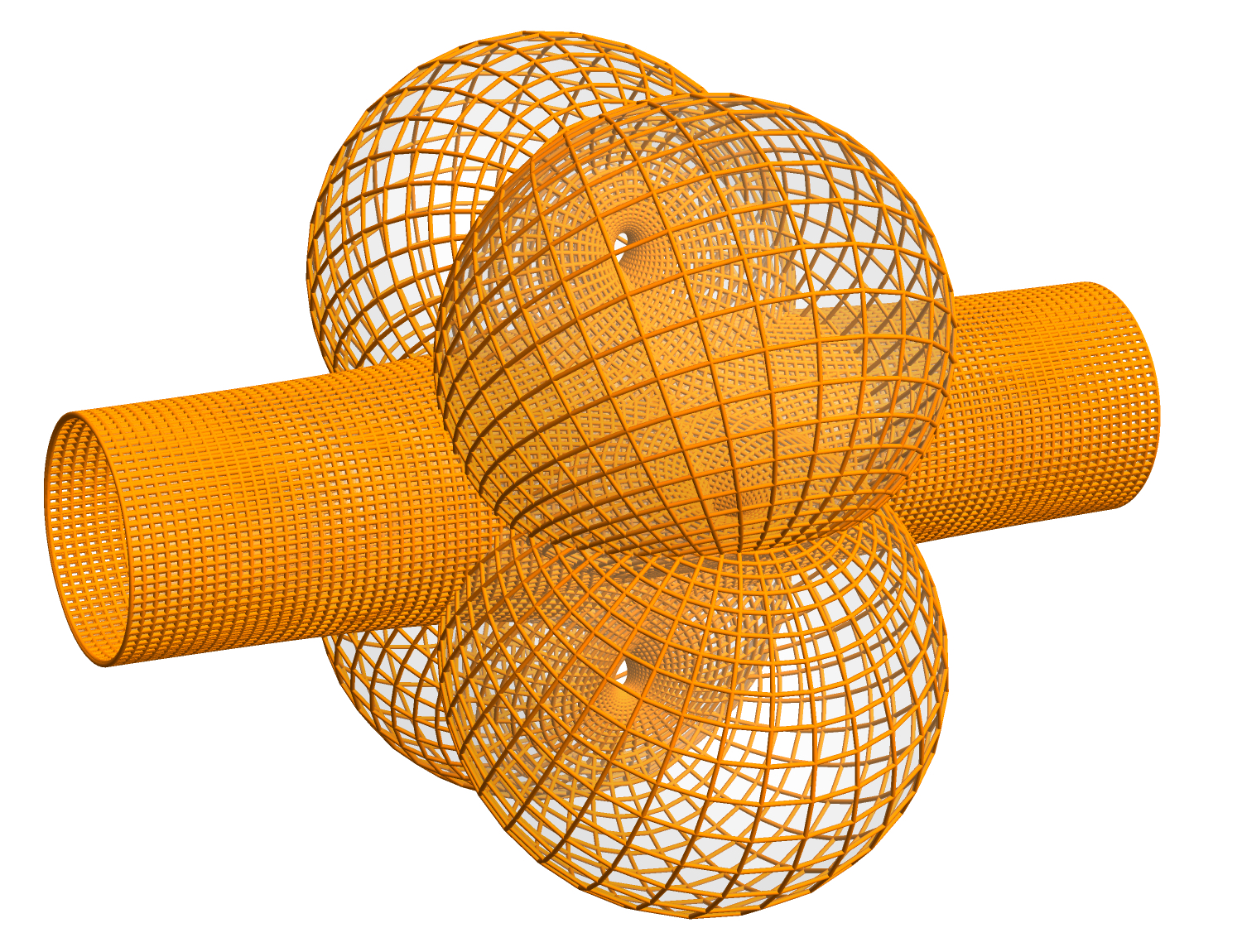}
	\end{minipage}
	\caption{Discrete cmc bubbletons drawn using explicit parametrisation \eqref{eqn:expPar} with initial condition as in \eqref{eqn:cmcInit} (on the left: $k = 3$, $\rho = 1$, $M = 120$, $N = 15$; on the right: $k = 4$, $\rho = 3$, $M = 80$, $N = 10$).}
	\label{fig:bubbleton2}
\end{figure}

\begin{remark}
	The condition that  $M$ is sufficiently large is essential, setting the discrete case apart from the smooth case.
	In the smooth case, the condition $k > \rho$ is sufficient to guarantee the existence of cmc cylinders appearing as Darboux transforms of the circular cylinder; however in the discrete case, there are configurations where $k > \rho$ is not sufficient.
	For example, if we choose $k = 5$, $\rho = 2$, and $M = 3$, then we have
	\[
		0 < \nu = \frac{1}{4}\left(1 - \cot^2 \frac{\pi}{3}\tan^2 \frac{5 \pi}{6} \right) = \frac{2}{9} < \frac{1}{4}.
	\]
\end{remark}

\begin{figure}
	\centering
	\begin{minipage}{0.46\linewidth}
		\centering
		\includegraphics[width=0.9\textwidth]{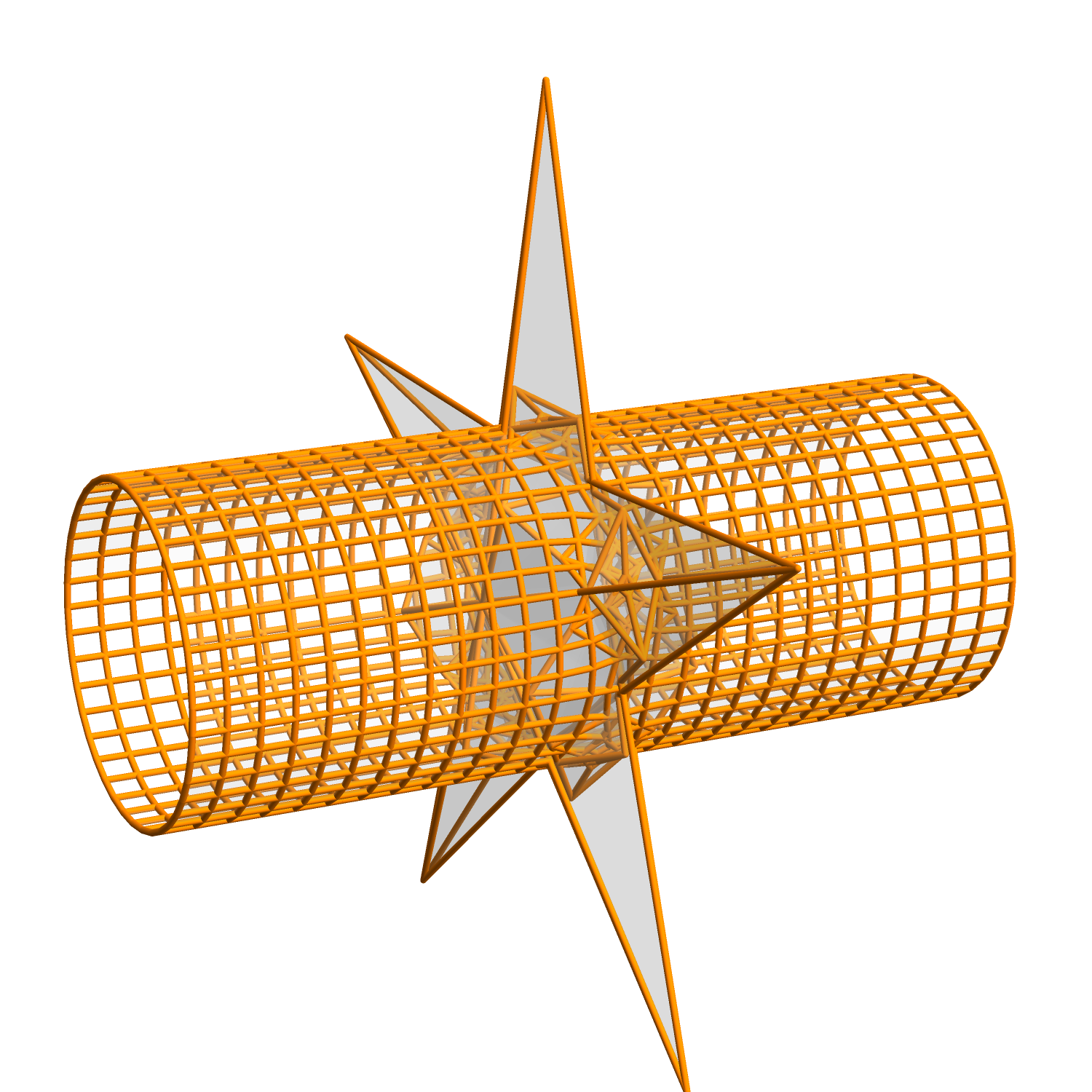}
	\end{minipage}
	\begin{minipage}{0.46\linewidth}
		\centering
		\includegraphics[width=0.9\textwidth]{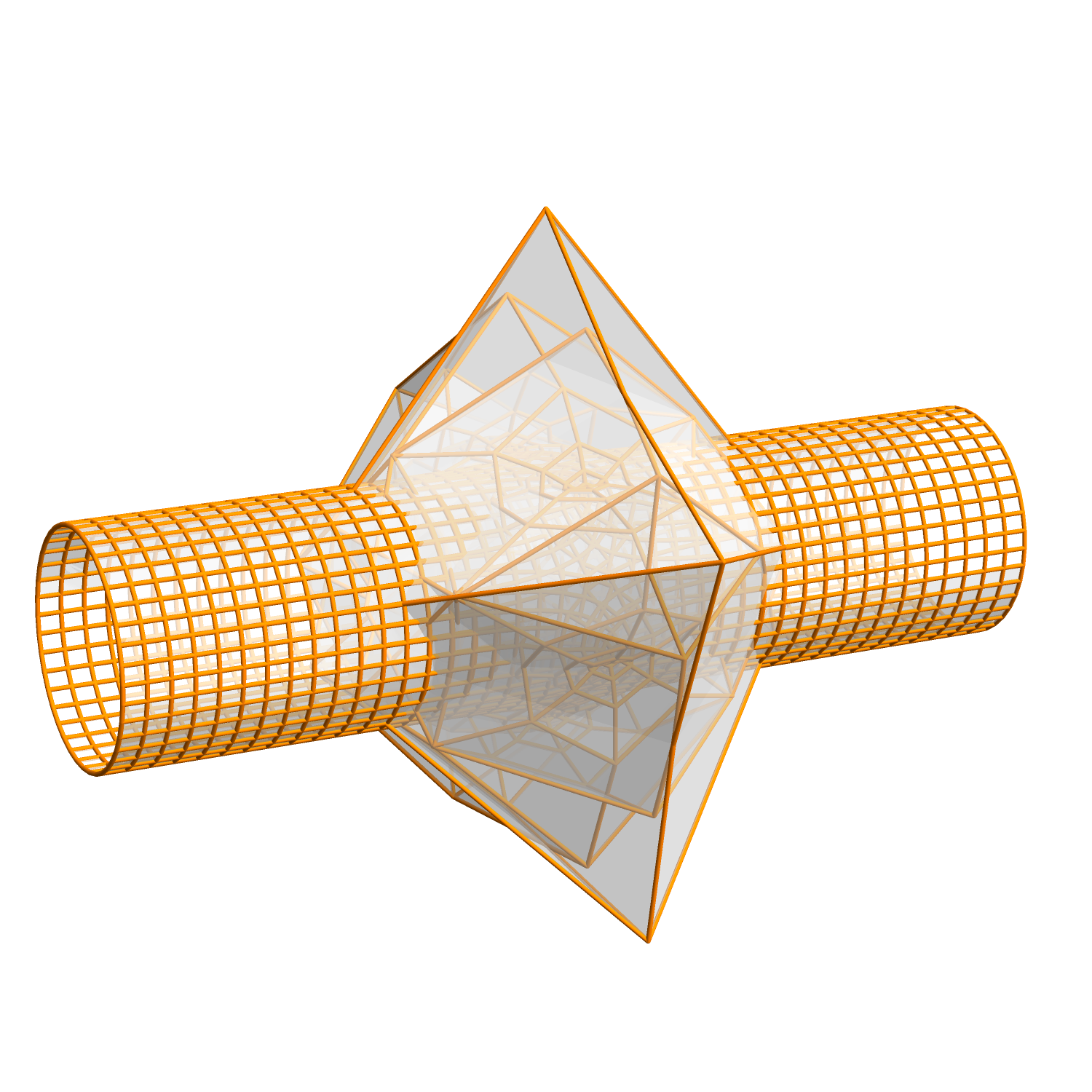}
	\end{minipage}
	\begin{minipage}{0.46\linewidth}
		\centering
		\includegraphics[width=0.9\textwidth]{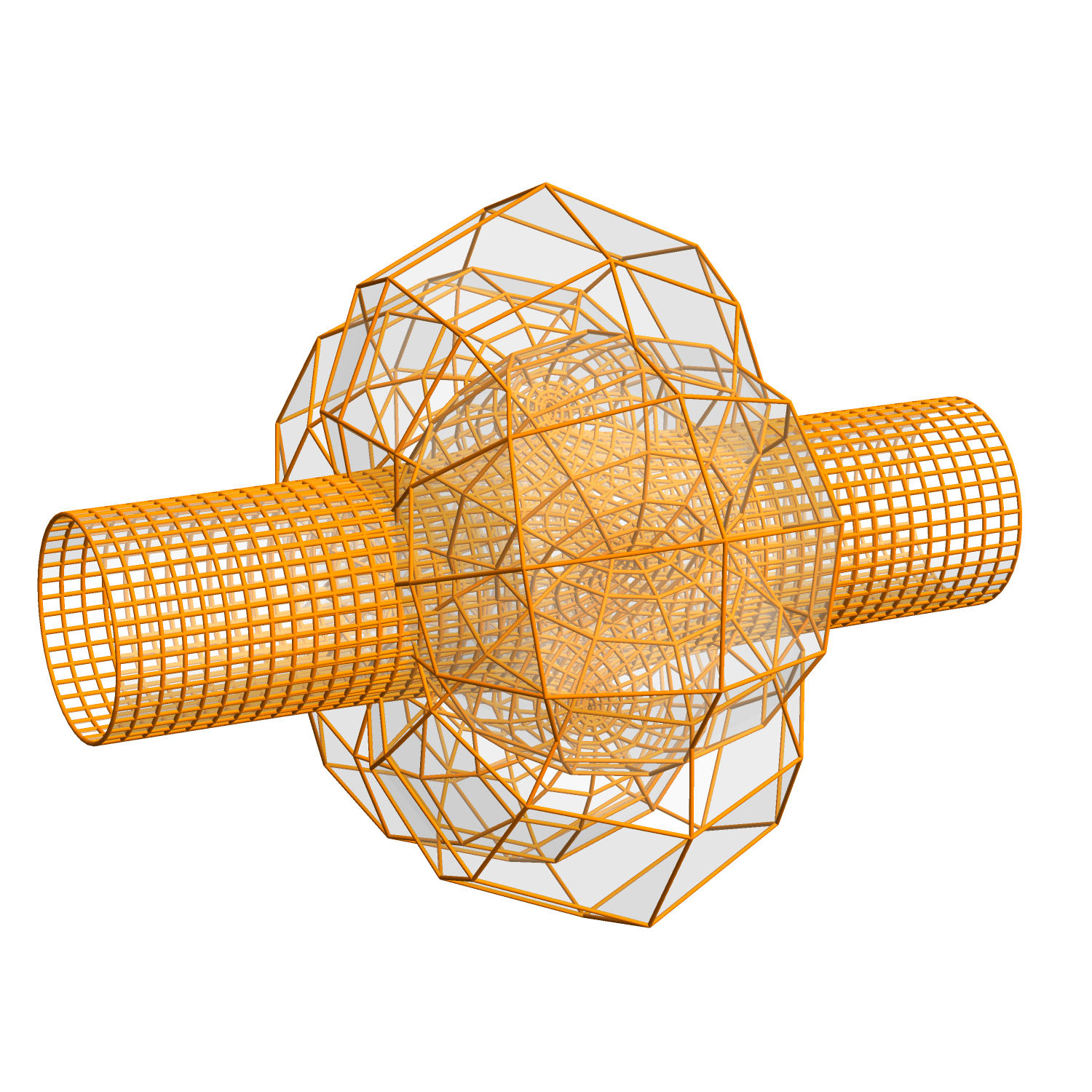}
	\end{minipage}
	\begin{minipage}{0.46\linewidth}
		\centering
		\includegraphics[width=0.9\textwidth]{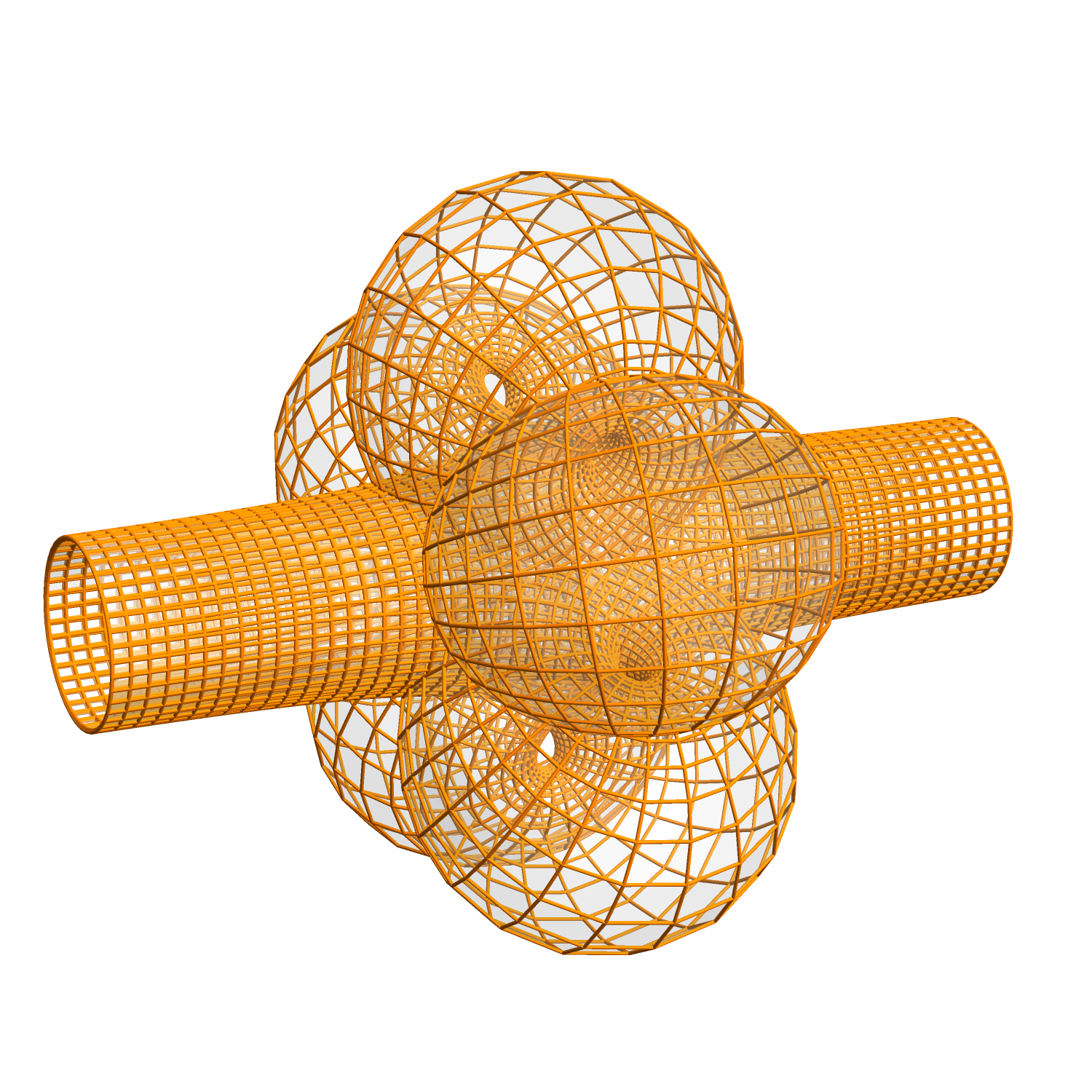}
	\end{minipage}
	\caption{Discrete cmc bubbletons drawn using explicit parametrisation \eqref{eqn:expPar} where $k = 5$, $M = 40$, and $N = 5$ with initial condition as in \eqref{eqn:cmcInit} at varying number of covers (top row from the left: $\rho = 1, 2$; bottom row from the left: $\rho = 3, 4$).}
	\label{fig:bubbleton3}
\end{figure}

To obtain closed-form parametrisations of discrete cmc bubbletons, the fact that $\nu < 0$ allows us to use Example~\ref{exam:explicit}; therefore, let $\hat{f}$ be parametrised by \eqref{eqn:expPar} so that via \eqref{eqn:yay}, we have
	\[
		\hat{f}_{0,0} = - \jj \, \frac{ \sin \frac{(\rho + k) \pi}{\rho M} c_2 + \sin \frac{(\rho - k) \pi}{\rho M}}{\sin \frac{(\rho - k) \pi}{\rho M} c_2 +  \sin \frac{(\rho + k) \pi}{\rho M}}.
	\]
For simplicity, we will assume that $c_2 \in \mathbb{R}$ so that $\hat{f}_{0,0} = \jj q$ for some $q \in \mathbb{R}$.
Thus, we need to find $c_2$ such that
	\[
		4 - \frac{1}{\nu} = |\hat{f}_{0,0} - f^*_{0,0}|^2 = \frac{4 (c_2 - 1)^2 \cos^2 \frac{\pi}{M}\sin^2 \frac{k\pi}{\rho M}}{\big(\sin \frac{(\rho - k) \pi}{\rho M} c_2 + \sin \frac{(\rho + k) \pi}{\rho M}\big)^2},
	\]
and a direct calculation shows
	\begin{equation}\label{eqn:cmcInit}
		c_2 = \pm \frac{\sqrt{\cos \frac{2\pi}{M} - \cos \frac{2 k \pi}{\rho M}}}{\sqrt{2} \sin \frac{(\rho - k) \pi}{\rho M}}.
	\end{equation}
For such $c_2$, we can use the closed-form parametrisation given in \eqref{eqn:expPar} to obtain discrete cmc bubbletons.
See Figure~\ref{fig:bubbleton}, \ref{fig:bubbleton2}, and \ref{fig:bubbleton3} for examples.

\subsection{Discrete isothermic tori}\label{sect:tori}
To construct non-trivial examples of discrete isothermic tori, let us consider the discrete homogeneous tori in $\mathbb{S}^3$ as our given surface $f$ so that
	\[
		f_{m,n} := q e^{\frac{2 \pi \ii}{N}n} + \jj p e^{\frac{2 \pi \ii}{M}m}
	\]
for non-zero $M,N \in\mathbb{Z}$ and $p,q \in \mathbb{R}$ satisfying $p^2 + q^2 = 1$.
We take the $\rho_1$--fold cover in the $m$-direction and $\rho_2$--fold cover in the $n$-direction.
Then we calculate that
	\[
		\cratio(f_i, f_j, f_k, f_\ell) = -\frac{p^2 |1 - e^\frac{2\pi \ii}{M}|^2}{q^2 |1 - e^\frac{2\pi \ii}{N}|^2} = \frac{\mu_{i\ell}}{\mu_{ij}}.
	\]
Thus we assume without loss of generality that
	\[
		\frac{1}{\mu_{ij}} = p^2 |1 - e^\frac{2\pi \ii}{M}|^2, \quad \frac{1}{\mu_{i\ell}} = -q^2 |1 - e^\frac{2\pi \ii}{N}|^2,
	\]
so that every $m$-curvature line is arc-length polarised while every $n$-curvature line is negative arc-length polarised.

Choosing any base point $(m_0, n_0)$, let us first take the $m$-curvature line $f_{m, n_0}$ through the base point.
Then since we have a discrete circle with constant polarisation, we deduce that any of its Darboux transforms $\hat{f}_{m, n_0}$ with spectral parameter
	\[
		\nu_1 = \frac{1}{4p^2}\left(1 - \cot^2 \frac{\pi}{M} \tan^2 \frac{k_1 \pi}{\rho_1 M}\right)
	\]
for any $k_1 \in \mathbb{Z}$ results in a closed discrete polarised curve.
Similarly, the $n$-curvature line $f_{m_0, n}$ is also a discrete circle with constant polarisation, so that any of its Darboux transforms $\hat{f}_{m_0, n}$ with spectral parameter
	\[
		\nu_2 = -\frac{1}{4q^2}\left(1 - \cot^2 \frac{\pi}{N} \tan^2 \frac{k_2 \pi}{\rho_2 N}\right)
	\]
for $k_2 \in \mathbb{Z}$ results in a closed discrete polarised curve.

Thus, if we choose $p$ and $q$ so that
	\begin{equation}\label{eqn:torires}
		\begin{aligned}
			p^2 &= \frac{1- \cot^2 \frac{\pi}{M} \tan^2 \frac{k_1 \pi}{\rho_1 M}}{\cot^2 \frac{\pi}{N} \tan^2 \frac{k_2 \pi}{\rho_2 N} - \cot^2 \frac{\pi}{M} \tan^2 \frac{k_1 \pi}{\rho_1 M}},\\
			q^2 &= \frac{1 - \cot^2 \frac{\pi}{N} \tan^2 \frac{k_2 \pi}{\rho_2 N}}{\cot^2 \frac{\pi}{M} \tan^2 \frac{k_1 \pi}{\rho_1 M} - \cot^2 \frac{\pi}{N} \tan^2 \frac{k_2 \pi}{\rho_2 N}},
		\end{aligned}
	\end{equation}
then we have $p^2 + q^2 = 1$, and
	\begin{equation}\label{eqn:torisp}
		\nu := \nu_1 = \nu_2 = \frac{1}{4}\left(\cot^2 \frac{\pi}{N} \tan^2 \frac{k_2 \pi}{\rho_2 M} - \cot^2 \frac{\pi}{M} \tan^2 \frac{k_1 \pi}{\rho_1 M}\right).
	\end{equation}
This implies that for such $p$ and $q$, we have that every Darboux transform with such spectral parameter $\nu = \nu_1 = \nu_2$ gives a discrete isothermic torus.
These are the discrete analogues of the isothermic tori found in \cite{bernstein_non-special_2001} (see also \cite{hertrich-jeromin_introduction_2003}).

To obtain closed-form parametrisations of the discrete isothermic tori, let us take $p$ and $q$ as in \eqref{eqn:torires} so that the spectral parameter satisfies
	\[
		\nu = \frac{1}{4p^2}\left(1 - \cot^2 \frac{\pi}{M} \tan^2 \frac{k_1 \pi}{\rho_1 M}\right) = -\frac{1}{4q^2}\left(1 - \cot^2 \frac{\pi}{N} \tan^2 \frac{k_2 \pi}{\rho_2 N}\right).
	\]
Since $f_{m,0} = q + \jj p e^{\frac{2\pi \ii}{M}}m$, we use the discussions from Section~\ref{sect:circle} to obtain that
	\begin{equation}\label{eqn:toriAB0}
		\begin{aligned}
			a_{m,0} &= \left(e^{-\frac{\pi \ii}{M}} \cos \tfrac{\pi}{M} \sec \tfrac{k_1 \pi}{\rho_1 M}\right)^m \alpha_m, \\
			b_{m,0} &= \tfrac{1}{2p}  \csc{\tfrac{\pi}{M}} \sec{\tfrac{k_1 \pi}{\rho_1 M}} \left(e^{-\frac{\pi \ii}{M}} \cos \tfrac{\pi}{M} \sec \tfrac{k_1 \pi}{\rho_1 M}\right)^m \beta_m
		\end{aligned}
	\end{equation}
where
	\begin{align*}
		\alpha_m &= e^{-\frac{k_1 \pi \ii}{\rho_1 M}m}c^+ + e^{\frac{k_1 \pi \ii}{\rho_1 M}m}c^-, \\
		\beta_m &= \jj \left(e^{-\frac{k_1 \pi \ii}{\rho_1 M}m} \sin{\tfrac{(\rho_1 + k_1) \pi}{\rho_1 M}} c^+ + e^{\frac{k_1 \pi_1 \ii}{\rho_1 M}m} \sin{\tfrac{(\rho_1 - k_1) \pi}{\rho_1 M}}  c^-\right)
	\end{align*}
for some $c^\pm \in \mathbb{H}$.

Now to solve along an $n$-curvature line direction, we fix some $m$ and note that $n$-curvature lines are also circles; thus, we can use a completely analogous argument in Section~\ref{sect:circle} to find that
	\begin{equation}\label{eqn:toriAB}
		\begin{aligned}
			a_{m,n} &= \big(e^{\frac{\pi \ii}{N}} \cos \tfrac{\pi}{N} \sec \tfrac{k_2 \pi}{\rho_2 N}\big)^n \Big(e^{\frac{k_2 \pi \ii}{\rho_2 N}n}\gamma^+_m + e^{-\frac{k_2 \pi \ii}{\rho_2 N}n}\gamma^-_m\Big), \\
			b_{m,n} &=
				\begin{multlined}[t]
					-\tfrac{1}{2q}  \csc{\tfrac{\pi}{N}} \sec{\tfrac{k_2 \pi}{\rho_2 N}} \big(e^{-\frac{\pi \ii}{N}} \cos \tfrac{\pi}{N} \sec \tfrac{k_2 \pi}{\rho_2 N}\big)^n \\
					\qquad\Big(e^{\frac{k_2 \pi \ii}{\rho_2 N}n}\sin\tfrac{(\rho_2 + k_2) \pi}{\rho_2 N} \gamma^+_m + e^{-\frac{k_2 \pi \ii}{\rho_2 N}n}\sin\tfrac{(\rho_2 - k_2) \pi}{\rho_2 N}\gamma^-_m\Big)
				\end{multlined}
		\end{aligned}
	\end{equation}
for some quaternion valued functions $\gamma^\pm_m$ of $m$.
To find $\gamma^\pm_m$, we let $n = 0$ in \eqref{eqn:toriAB} to see that
	\begin{align*}
		a_{m,0} &= \gamma^+_m + \gamma^-_m, \\
		b_{m,0} &= -\tfrac{1}{2q}  \csc{\tfrac{\pi}{N}} \sec{\tfrac{k_2 \pi}{\rho_2 N}} \Big(\sin\tfrac{(\rho_2 + k_2) \pi}{\rho_2 N} \gamma^+_m + \sin\tfrac{(\rho_2 - k_2) \pi}{\rho_2 N}\gamma^-_m\Big).
	\end{align*}
Solving for $\gamma^\pm_m$ and using \eqref{eqn:toriAB0}, we obtain that
	\begin{align*}
		\gamma^\pm_m
			&= \!\begin{multlined}[t]
					\mp \tfrac{1}{2}\sec{\tfrac{\pi}{N}} \csc{\tfrac{k_2 \pi}{\rho_2 N}} \big(e^{-\frac{\pi \ii}{M}} \cos \tfrac{\pi}{M} \sec \tfrac{k_1 \pi}{\rho_1 M}\big)^m \\
					\qquad \Big( \sin\tfrac{(\rho_2 \mp k_2) \pi}{\rho_2 N} \alpha_m + \tfrac{q}{p} \sin{\tfrac{\pi}{N}} \cos{\tfrac{k_2 \pi}{\rho_2 N}} \csc{\tfrac{\pi}{M}} \sec{\tfrac{k_1 \pi}{\rho_1 N}} \beta_m \Big)
			\end{multlined}\\
			&=: \mp \tfrac{1}{2} \sec{\tfrac{\pi}{N}} \csc{\tfrac{k_2 \pi}{\rho_2 N}}  \big(e^{-\frac{\pi \ii}{M}} \cos \tfrac{\pi}{M} \sec \tfrac{k_1 \pi}{\rho_1 M}\big)^m \tilde{\gamma}^\pm_m. 
	\end{align*}
Therefore,
	\begin{align*}
		a_{m,n} &=
			\!\begin{multlined}[t]
				-\tfrac{1}{2}\sec{\tfrac{\pi}{N}} \csc{\tfrac{k_2 \pi}{\rho_2 N}} \big(e^{-\frac{\pi \ii}{M}} \cos \tfrac{\pi}{M} \sec \tfrac{k_1 \pi}{\rho_1 M}\big)^m \big(e^{\frac{\pi \ii}{N}} \cos \tfrac{\pi}{N} \sec \tfrac{k_2 \pi}{\rho_2 N}\big)^n \\
				\qquad \Big(e^{\frac{k_2 \pi \ii}{\rho_2 N}n}\tilde{\gamma}^+_m - e^{-\frac{k_2 \pi \ii}{\rho_2 N}n}\tilde{\gamma}^-_m\Big)
			\end{multlined}\\
		b_{m,n} &=
			\!\begin{multlined}[t]
				\tfrac{1}{q}\csc \tfrac{2 \pi}{N}\csc \tfrac{2 k_2 \pi}{\rho_2 N} \big(e^{-\frac{\pi \ii}{M}} \cos \tfrac{\pi}{M} \sec \tfrac{k_1 \pi}{\rho_1 M}\big)^m \big(e^{-\frac{\pi \ii}{N}} \cos \tfrac{\pi}{N} \sec \tfrac{k_2 \pi}{\rho_2 N}\big)^n \\
				\qquad\Big(e^{\frac{k_2 \pi \ii}{\rho_2 N}n}\sin\tfrac{(\rho_2 + k_2) \pi}{\rho_2 N}\tilde{\gamma}^+_m - e^{-\frac{k_2 \pi \ii}{\rho_2 N}n}\sin\tfrac{(\rho_2 - k_2) \pi}{\rho_2 N}\tilde{\gamma}^-_m\Big).
 			\end{multlined}
	\end{align*}
Summarzing, we have
	\begin{equation}\label{eqn:abinvtori}
		a_{m,n} b_{m,n}^{-1} =
			\!\begin{multlined}[t]
				-2q \sin \tfrac{\pi}{N} \cos \tfrac{k_2 \pi}{\rho_2 N} e^{-\frac{\pi \ii}{M}m}e^{\frac{\pi \ii}{N}n} \Big(e^{\frac{k_2 \pi \ii}{\rho_2 N}n}\tilde{\gamma}^+_m - e^{-\frac{k_2 \pi \ii}{\rho_2 N}n}\tilde{\gamma}^-_m\Big) \\
				\Big(e^{\frac{k_2 \pi \ii}{\rho_2 N}n}\sin\tfrac{(\rho_2 + k_2) \pi}{\rho_2 N}\tilde{\gamma}^+_m - e^{-\frac{k_2 \pi \ii}{\rho_2 N}n}\sin\tfrac{(\rho_2 - k_2) \pi}{\rho_2 N}\tilde{\gamma}^-_m\Big)^{-1}e^{\frac{\pi \ii}{M}m}e^{\frac{\pi \ii}{N}n},
			\end{multlined}
	\end{equation}
where
	\[
		\tilde{\gamma}^\pm_m = \sin\tfrac{(\rho_2 \mp k_2) \pi}{\rho_2 N} \alpha_m + \tfrac{q}{p} \sin{\tfrac{\pi}{N}} \cos{\tfrac{k_2 \pi}{\rho_2 N}} \csc{\tfrac{\pi}{M}} \sec{\tfrac{k_1 \pi}{\rho_1 M}} \beta_m
	\]
for 
	\begin{align*}
		\alpha_m &= e^{-\frac{k_1 \pi \ii}{\rho_1 M}m}c^+ + e^{\frac{k_1 \pi \ii}{\rho_1 M}m}c^-, \\
		\beta_m &= \jj \left(e^{-\frac{k_1 \pi \ii}{\rho_1 M}m} \sin{\tfrac{(\rho_1 + k_1) \pi}{\rho_1 M}} c^+ + e^{\frac{k_1 \pi_1 \ii}{\rho_1 M}m} \sin{\tfrac{(\rho_1 - k_1) \pi}{\rho_1 M}}  c^-\right),
	\end{align*}
and the discrete isothermic tori $\hat{f}$ obtained via Darboux transformations of the homogeneous tori are given by
	\[
		\hat{f}_{m,n} = f_{m,n} + a_{m,n} b_{m,n}^{-1}.
	\]

If the constants $c^\pm \in \mathbb{H}$ are chosen so that for any initial point $(m_0, n_0)$, we have $\hat{f}_{m_0, n_0} \in \mathbb{S}^3$, then it follows that $\hat{f}_{m,n}$ is a discrete isothermic torus in $\mathbb{S}^3$ (see Fact~\ref{fact:s3}).

\begin{figure}
	\centering
	\begin{minipage}{0.46\linewidth}
		\centering
		\includegraphics[width=0.9\textwidth]{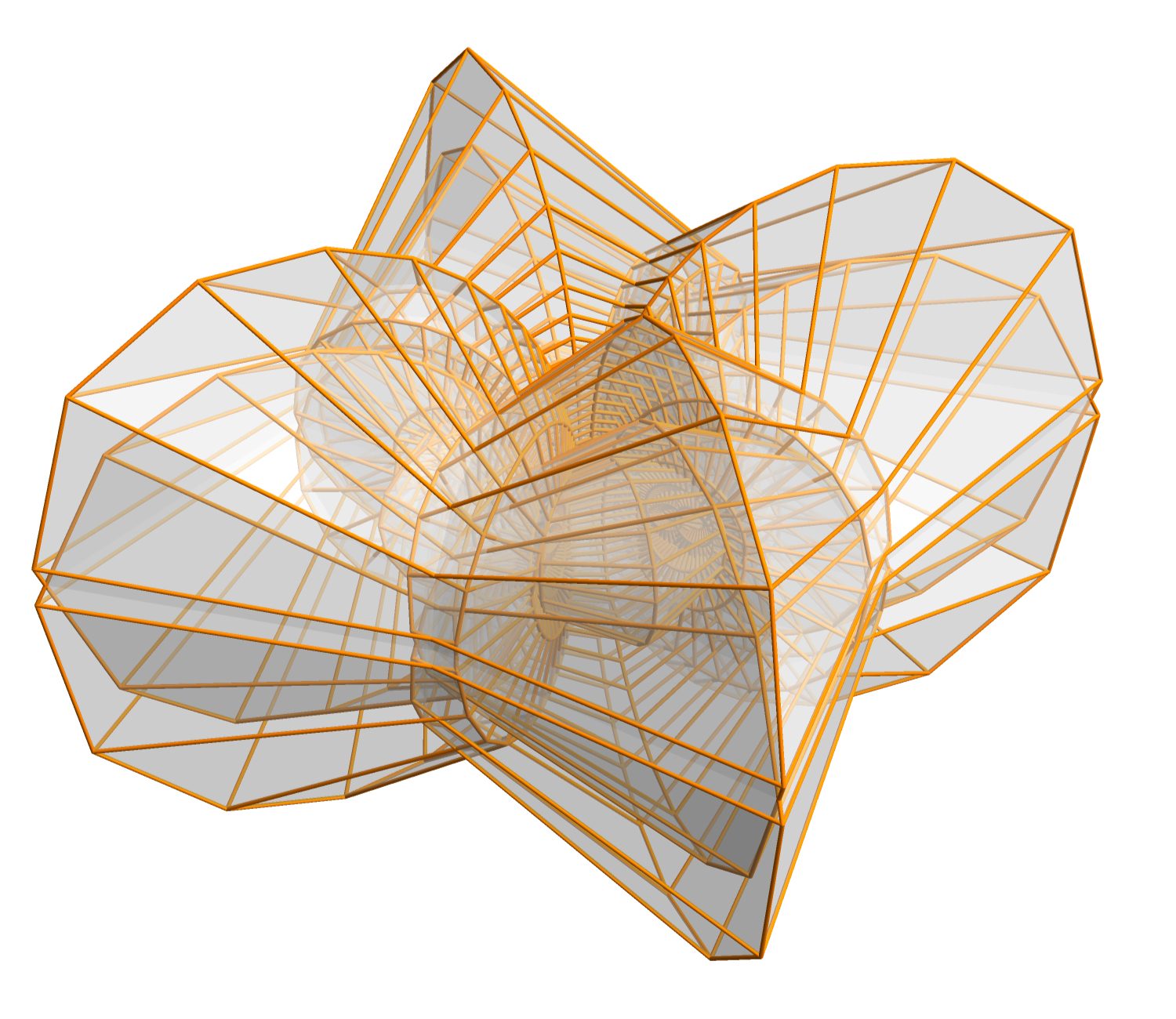}
	\end{minipage}
	\begin{minipage}{0.46\linewidth}
		\centering
		\includegraphics[width=0.9\textwidth]{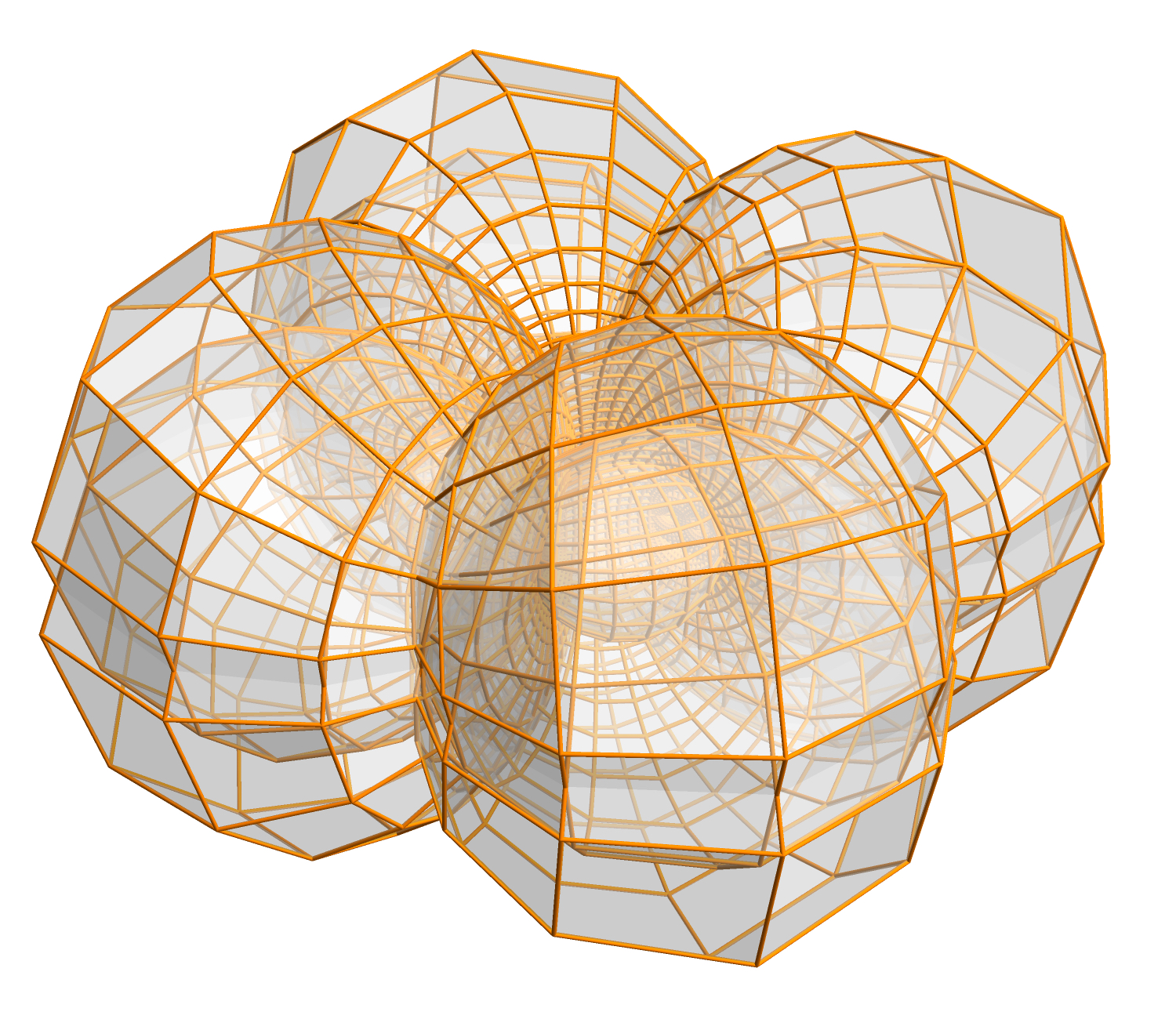}
	\end{minipage}
	\begin{minipage}{0.46\linewidth}
		\centering
		\includegraphics[width=0.9\textwidth]{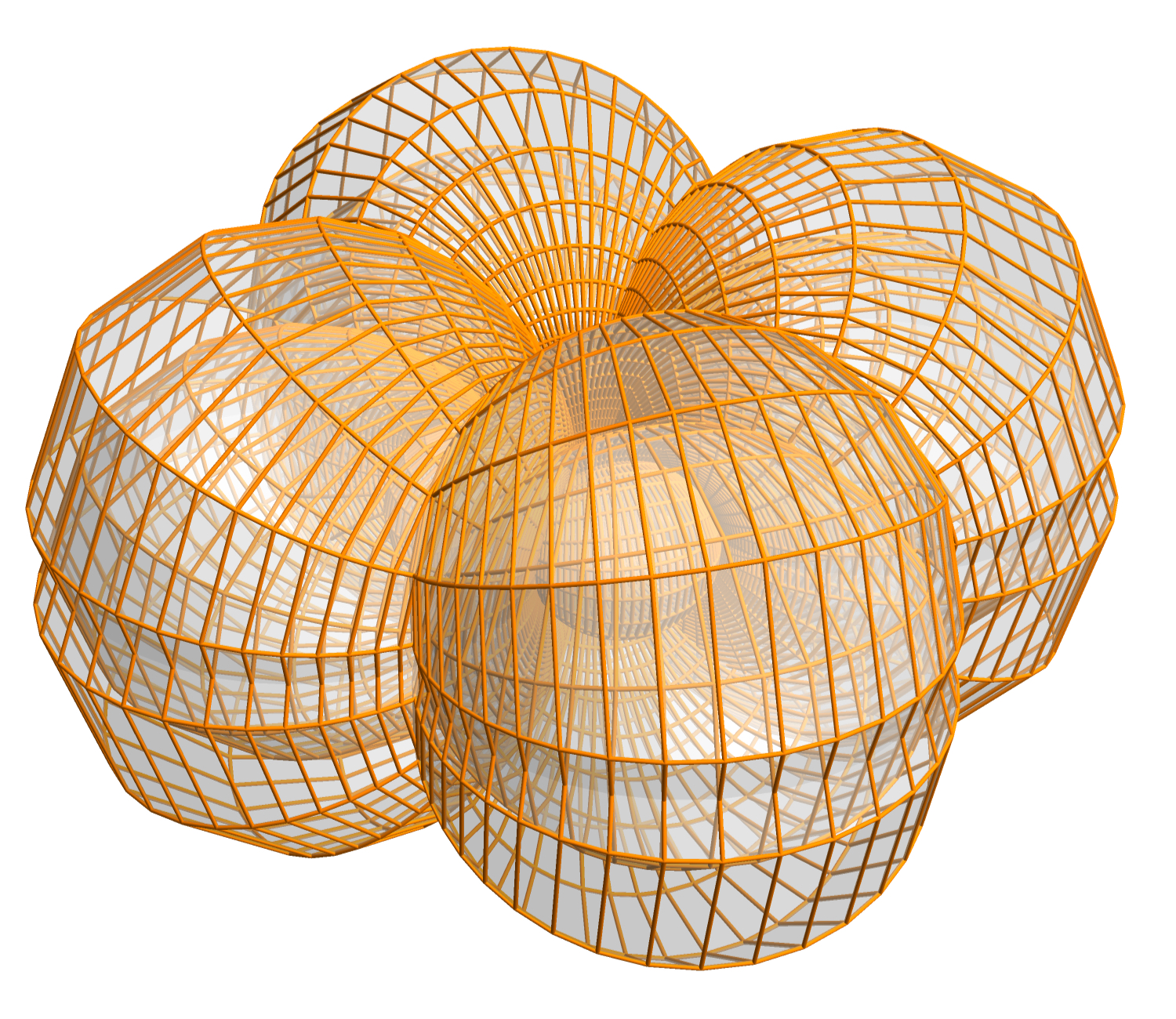}
	\end{minipage}
	\begin{minipage}{0.46\linewidth}
		\centering
		\includegraphics[width=0.9\textwidth]{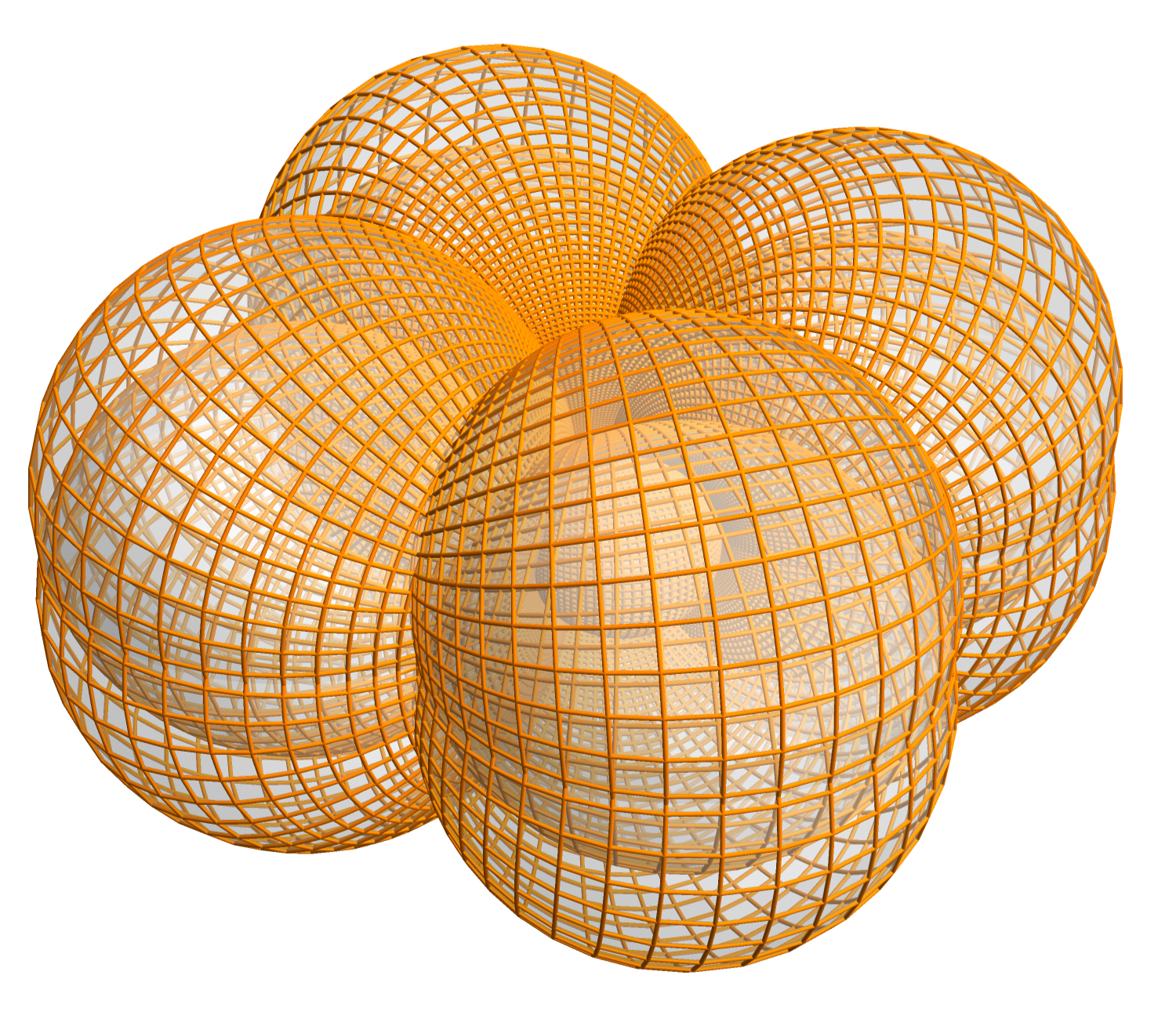}
	\end{minipage}
	\caption{Discrete isothermic tori in $\mathbb{S}^3$ drawn using the explicit parametrisation \eqref{eqn:abinvtori} with $k_1 = 4$, $\rho_1 = 3$, $k_2 = 2$, $\rho_2 = 3$, $c^+ = 0.45 + \jj r_2$, and $c^- = 1$ with different number of subdivisions and corresponding values of $r_2$ (on the top left: $M = 12$, $N = 40$, $r_2 \approx 0.395155$; on the top right: $M = 40$, $N = 40$, $r_2 \approx 0.385119$; on the bottom left: $M = 160$, $N = 40$, $r_2 \approx 0.384212$; on the bottom right: $M = 160$, $N = 160$, $r_2 \approx 0.384438$).}
	\label{fig:tori}
\end{figure}

\begin{figure}
	\centering
	\begin{minipage}{0.46\linewidth}
		\centering
		\includegraphics[width=0.8\textwidth]{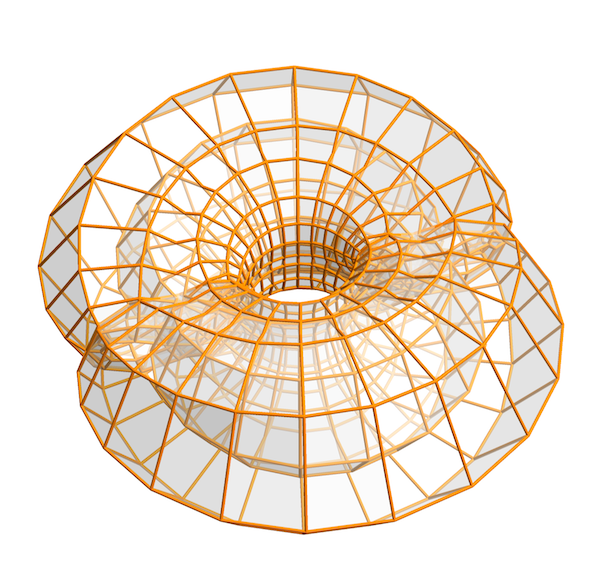}
	\end{minipage}
	\begin{minipage}{0.46\linewidth}
		\centering
		\includegraphics[width=\textwidth]{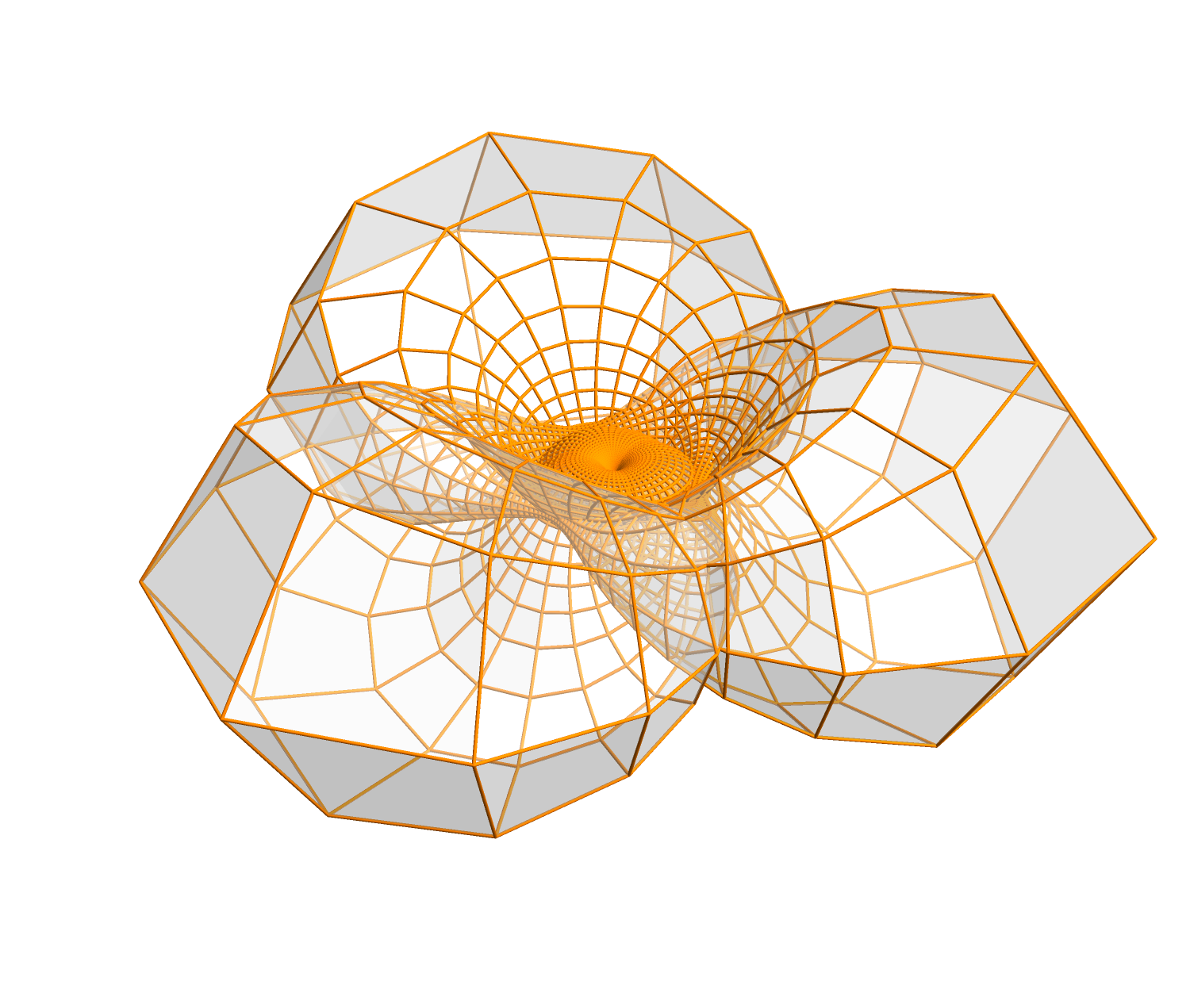}
	\end{minipage}
	\caption{Discrete isothermic tori in $\mathbb{S}^3$ drawn using the explicit parametrisation \eqref{eqn:abinvtori}
	(on the left: $k_1 = 2$, $\rho_1 = 3$, $k_2 = 3$, $\rho_2 = 2$, $M = 12$, $N = 12$, $c^+ \approx 0.45 + \jj 0.236786$, $c^- = 1$; on the right: $k_1 = 3$, $\rho_1 = 2$, $k_2 = 2$, $\rho_2 = 3$, $M = 60$, $N = 40$, $c^+ \approx - 0.6 - \jj 1.70678$, $c^- = 1$).}
	\label{fig:tori2}
\end{figure}

\subsection{Continuum limit of the resonance points}\label{sec:continuum}
Note that the calculation of resonance points in each section is dependent on the number of subdivisions of the smooth circle.
Therefore, for the resonance points, taking the continuum limit corresponds to taking the limit as the number of subdivisions tend to infinity.
As we have investigated the monodromy problem of the discrete isothermic surfaces under structure-preserving discretisation, it is fair to expect that the resonance points of the discrete objects also converge to those of the smooth counterpart under the continuum limit.
We can check the convergence explicitly as follows.

For the surface of revolution case, we see that
	\[
		\lim_{M \to \infty} \frac{1}{4}\left(1 - \cot^2 \frac{\pi}{M} \tan^2 \frac{k \pi}{\rho M}\right) = \frac{\rho^2 - k^2}{4 \rho^2}
	\]
for $k\in \mathbb{Z}$, which is identical to the resonance point of surface of revolutions in the smooth case (see, for example, \cite[Section~2]{cho_generalised_2022}).

Also for the homogeneous tori case, we can similarly check that resonance points of the discrete case \eqref{eqn:torisp} converge to that of the smooth case:
	\begin{align*}
		\lim_{M,N \to \infty} \nu
			&= \lim_{M,N \to \infty} \frac{1}{4}\left(\cot^2 \frac{\pi}{N} \tan^2 \frac{k_2 \pi}{\rho_2 N} - \cot^2 \frac{\pi}{M} \tan^2 \frac{k_1 \pi}{\rho_1 M}\right) \\
			&= \frac{1}{4}\left(\frac{k_2^2}{\rho_2^2} - \frac{k_1^2}{\rho_1^2}\right).
	\end{align*}
In fact, it follows that even the radii of the discrete homogeneous tori admitting isothermic tori as Darboux transforms converge to the radii of the smooth case.
We can check this explicitly by noting that
	\[
		\lim_{M,N \to \infty}(1 - 4 \nu p^2) = \frac{k_1^2}{\rho_1^2}, \quad \lim_{M,N \to \infty} (1 + 4 \nu q^2) = \frac{k_2^2}{\rho_2^2},
	\]
which is consistent with the smooth case in \cite[Equation~(5.24)]{hertrich-jeromin_introduction_2003}.

These observations suggest that the monodromy problem of Darboux transforms of smooth isothermic surfaces can be investigated via structure preserving discretisation without any numerical integration.



\end{document}